\documentclass{amsart}

\usepackage{style}

\title{Counting elliptic fibrations on K3 surfaces}

\author{Dino Festi}
\address[Dino Festi]
{Dipartimento di matematica ``Federigo Enriques'', Università degli Studi di Milano, via Saldini~50, 20133 Milan, Italy}
\email{dino.festi@unimi.it}

\author{Davide Cesare Veniani}
\address[Davide Cesare Veniani]
{Institut für Geometrie und Topologie, 
Universität Stuttgart,
Pfaffenwaldring~57,
70569 Stuttgart, Germany}
\email{davide.veniani@mathematik.uni-stuttgart.de}

\date{\today}

\makeatletter
\@namedef{subjclassname@2020}{%
 \textup{2020} Mathematics Subject Classification}
\makeatother



\begin{document}

\maketitle

\begin{abstract}
We solve the problem of counting jacobian elliptic fibrations on an arbitrary complex projective K3 surface up to automorphisms.
We then illustrate our method with several explicit examples.
\end{abstract}

\section{Introduction}

An \emph{elliptic pencil} on a complex projective K3 surface \(X\) is a complete linear system of divisors whose general member is a smooth elliptic curve \(E\). 
An elliptic pencil corresponds to an \emph{elliptic fibration} \(\pi \colon X \rightarrow \IP^1\), which is said to be \emph{jacobian} if it admits a section \(O\). 
We denote by \(\cJ_X\) the set of all jacobian elliptic fibrations on \(X\).
The automorphism group \(\Aut(X)\) acts on \(\cJ_X\) with a finite number of orbits according to a result by Sterk \cite[Cor.~2.7]{Sterk:finiteness.results.K3}.
The aim of this paper is to determine the number of orbits
\[
    |\cJ_X/{\Aut(X)}|.
\]

Let \(\pi\colon X \to \IP^1\) be a fixed jacobian filtration and denote by \(E\) its generic fiber.
Once an arbitrary section \(O\) is chosen as neutral element, the set of sections of \(\pi\) acquires the structure of a group called the \emph{Mordell--Weil group} of the fibration \(\pi\).
The classes of \(E\) and \(O\) in the Néron--Severi lattice~\(S_X\) induce an embedding \(\iota\colon \bU \hookrightarrow S_X\), where \(\bU\) denotes the hyperbolic unimodular even lattice of rank~\(2\).
The orthogonal complement 
\[
    W \coloneqq \iota(\bU)^\perp \subset S_X,
\]
is called the \emph{frame} of the fibration \(\pi\).
The frame~\(W\) encodes virtually all geometrical information about the fibration~\(\pi\).
Indeed, the dual graph of the components of the reducible fibers can be inferred from its root sublattice~\(W_\rootlattice\), and the Mordell--Weil group is isomorphic to \(W/W_\rootlattice\) (see for instance the survey by Schütt and Shioda~\cite{Schuett.Shioda:elliptic.surfaces}).
All frames have the same signature and discriminant form, so they belong to the same lattice genus \(\cW_X\), which we call the \emph{frame genus} of~\(X\) (\autoref{def:frame_genus}). 

Although several works have been dedicated to the classification of jacobian elliptic fibrations on a given K3 surface \cite{Bertin.Garbagnati.Hortsch.Lecacheux.Mase.Salgado.Whitcher, Bertin.Lecacheux:apery-fermi.2-isogenies, Braun.Kimura.Watari:classif.ell.fibr, Comparin.Garbagnati:vanGeemen-Sarti.inv, Elkies.Schuett:genus.1.fibrations, Harrache.Lecacheux:etude, Kloosterman:classification.jac.ell.fibr, Kumar:elliptic.fibrations, Mezzedimi:K3.zero.entropy, Nishiyama:Jacobian.fibrations.K3.MW, Oguiso:jacobian.fibrations.Kummer, Shimada:on.ell.K3},
in very few of them the number \(|\cJ_X/{\Aut(X)}|\) has been explicitly determined. 
The various approaches were clarified in an unpublished paper by Braun, Kimura and Watari~\cite{Braun.Kimura.Watari:classif.ell.fibr}, where the following distinctions were introduced:
\begin{itemize}
    \item the ``\(\cJ^{\mathrm{(type)}}(X)\) classification'' is the problem of determining the pairs \((W_\rootlattice,W/W_\rootlattice)\) for \(W \in \cW_X\);
    \item the ``\(\cJ_2(X)\) classification'' is the problem of determining the frame genus \(\cW_X\);
    \item the ``\(\cJ_1(X)\) classification'' is the problem of determining the orbits \(\cJ_X/{\Aut(X)}\).
    \item the ``\(\cJ_0(X)\) classification'' is the problem of determining \(\cJ_X\).
\end{itemize}
In general, each classification is coarser than the next one. 
For most geometrical applications, the ``\(\cJ^{\mathrm{(type)}}(X)\) classification'' suffices. 
Very few authors went as far as the ``\(\cJ_1(X)\) classification''. As for the ``\(\cJ_0(X)\) classification'', Nikulin \cite[Theorem 5.1]{Nikulin:ell.fibr.K3} proved that the set \(\cJ_X\) is finite if and only if \(S_X\) belongs to a certain finite set of lattices.

In a remarkable exception dating back to more than 30 years ago, Oguiso \cite{Oguiso:jacobian.fibrations.Kummer} proved that \[|\cJ_X/{\Aut(X)}| \in \set{16,23,38,59}\] when \(X\) is the Kummer surface associated to the product \(E \times F\) of two non-isogenous elliptic curves or, equivalently, when it holds \(T_X \cong \bU(2)^2\).
Oguiso's arguments are based on a deep understanding of the geometry of such Kummer surfaces and cannot be generalized to other K3 surfaces.

In the same unpublished paper~\cite{Braun.Kimura.Watari:classif.ell.fibr}, Braun, Kimura and Watari also found a sufficient condition for 
\(|\cJ_X/{\Aut(X)}| = |\cW_X|\) to hold (which can be seen as a direct corollary of our main theorem, see \autoref{cor:uniform.bound}), and listed a few cases where this condition is true.
One of these cases, namely the singular K3 surface \(X\) with transcendental lattice \(T_X \cong [2]\oplus[6]\), was later completely worked out by Bertin et al.~\cite{Bertin.Garbagnati.Hortsch.Lecacheux.Mase.Salgado.Whitcher, Bertin.Lecacheux:automorphisms.certain}, who found that
\[|\cJ_X/{\Aut(X)}| = |\cW_X| = 53.\]

The same sufficient condition was independently found by Mezzedimi~\cite{Mezzedimi:K3.zero.entropy}, who also proved that a K3 surface with Néron--Severi lattice \(S_X \cong \bU \oplus [-2d]\) satisfies \(|\cJ_X/{\Aut(X)}| = 1\) if \(d =1\) and
\[ 
    |\cJ_X/{\Aut(X)}| = 2^{k-1},
\] 
if \(d \geq 2\), where \(k\) is the number of prime divisors of \(d\) \cite[Proposition~4.4]{Mezzedimi:K3.zero.entropy}.
Moreover, he showed that \(|\cJ_X| = 1\) if \(S_X\) belongs to a certain explicit list of \(32\) lattices~\cite[Thm.~6.14]{Mezzedimi:K3.zero.entropy}.

We are not aware of other cases in the literature where \(|\cJ_X/{\Aut(X)}|\) is explicitly known.


Our main result, namely \autoref{thm:main_multiplicity}, is a formula for the number of jacobian fibrations up to automorphisms with the same frame \(W \in \cW_X\), which we call the \emph{multiplicity} of the frame~\(W\). 
The multiplicity of \(W\) turns out to be equal to the number of certain double cosets 
\[
    |H\backslash G /K|
\] 
in the orthogonal group \(G = \OO(T_X^\sharp)\), where \(T_X^\sharp\) is the discriminant group of the transcendental lattice \(T_X\). 
The subgroup \(H\) is related to the Hodge isometries of~\(T_X\), while the subgroup \(K\) depends on \(W\).
Quite interestingly, this pattern is shared by many other enumerative problems such as counting Kummer structures on~\(X\) (cf.~\cite{Hosono.Lian.Oguiso.Yau:Kummer.structures}), counting Fourier--Mukai partners of~\(X\) (cf.~\cite{Hosono.Lian.Oguiso.Yau:Fourier-Mukai.number.K3.surf}), or counting Enriques surfaces covered by~\(X\) (cf.~\cite{Shimada.Veniani:Enriques.involutions.singular.K3}).

The advantage of our algebraic method is that it can be implemented as soon as the transcendental lattice \(T_X\) is known. 
%
We compute \(|\cJ_X/{\Aut(X)}|\) explicitly in the following cases: 
\begin{itemize}
    \item K3 surfaces belonging to the Barth--Peters family (see \autoref{thm:Barth-Peters});
    \item Kummer surfaces associated to the product of non-isogenous elliptic curves, confirming Oguiso's results \cite{Oguiso:jacobian.fibrations.Kummer} (see \autoref{thm:Oguiso});
    \item Kummer surfaces associated to Jacobian of a very general curve of genus~\(2\), refining a work by Kumar~\cite{Kumar:elliptic.fibrations} (see \autoref{thm:Kumar});
    \item generic double covers of \(\IP^2\) ramified over \(6\) lines, refining a work by Kloosterman~\cite{Kloosterman:classification.jac.ell.fibr} (see \autoref{thm:Kloosterman});
    \item K3 surfaces belonging to the Apéry--Fermi pencil, refining a work by Bertin and Lecacheux~\cite{Bertin.Lecacheux:apery-fermi.2-isogenies} (see \autoref{thm:Apery-Fermi}).
\end{itemize}

In each case, all Gram matrices of the lattices in the frame genus are contained in the respective arXiv ancillary file.
Computations were carried out with GAP~\cite{GAP4}, Magma~\cite{Magma} and Sage~\cite{sage}.

\subsection*{Contents of the paper}
The paper is divided into two sections. 

The theoretical part is contained in \autoref{sec:main_thm}: it comprises the proof of our main theorem and general guidelines on how to implement our method.

All explicit examples are contained in \autoref{sec:examples}.

\subsection*{Acknowledgments} 
We warmly thank not only Simon Brandhorst for his help with Sage, but also Fabio Bernasconi, Alice Garbagnati and Remke Kloosterman for their insightful comments and Edgar Costa, Noam Elkies, Klaus Hulek, Giacomo Mezzedimi, Bartosz Naskr\k{e}cki, Matthias Schütt and Evgeny Shinder for their useful remarks on an earlier draft of this paper. 
The second author is also indebted to Christina Lienstromberg for her kind hospitality in Bonn, Germany, where this paper was finished.

\section{Main theorem} \label{sec:main_thm}

Throughout this section we let \(X\) be a complex projective K3 surface with Néron--Severi lattice~\(S_X\) and transcendental lattice~\(T_X\).

After fixing notation and conventions on lattices in~\autoref{sec:lattices} and recalling or proving some preliminary results in~\autoref{sec:preliminary_results}, we state the main theorem of the paper, namely \autoref{thm:main_multiplicity}, together with its immediate corollaries in \autoref{sec:statement}.
A proof of the theorem is given in~\autoref{sec:thm:main_mutliplicity}.
Finally, in~\autoref{sec:guidelines} we provide general guidelines on how to compute \(|\cJ_X/{\Aut(X)}|\) in explicit cases applying the formula of the theorem.

\subsection{Lattices} \label{sec:lattices}
A \emph{lattice} \(L\) of rank \(r\) is a free, finitely generated \(\IZ\)-module \(L \cong \IZ^r\) endowed with a symmetric bilinear form 
\[
L \times L \rightarrow \IZ, \quad (v, w) \mapsto v \cdot w\; .
\]

The \emph{dual} of \(L\) is the set
\[
    L^\vee \coloneqq \Set{ x \in L\otimes \IQ}{ x \cdot v \in \IZ \text{ for all \(v \in L\)}}
\]
and the \emph{discriminant group} is the finite abelian group 
\[
    L^\sharp \coloneqq L^\vee/L.
\]

If the lattice \(L\) is \emph{even}, meaning that \(v^2 \coloneqq v \cdot v \in 2\IZ\) for all \(v \in L\), then the form on \(L\) induces a finite quadratic form~\(L^\sharp \rightarrow \IQ/2\IZ\).
We write \(\OO(L)\) and \(\OO(L^\sharp)\) for the groups of isomorphisms of \(L\) and \(L^\sharp\) respecting the corresponding bilinear or quadratic forms. 
There is a natural homomorphism
\[
    \OO(L) \rightarrow \OO(L^\sharp), \quad \gamma \mapsto \gamma^\sharp,
\]
whose image is denoted \(\OO^\sharp(L)\).

An embedding \(\iota \colon M \hookrightarrow L\) of lattices is called \emph{primitive} if \(L/\iota(M)\) is a free group. A vector \(v \in L\) is called \emph{primitive} if \(\IZ v \hookrightarrow L\) is a primitive embedding.

A vector \(v\) in an even negative definite lattice \(L\) is called a \emph{root} if \(v^2 = -2\). The set of roots is denoted~\(\Delta(L)\). The sublattice generated by all roots is denoted \(L_\rootlattice\). 
A root \(v\) induces a reflection \(\rho_v \in \OO(L)\) defined by
\[
    \rho_v(w) \coloneqq w + (v\cdot w)v.
\]
The subgroup of \(\OO(L)\) generated by all reflections \(\rho_v\), denoted~\(\Weyl(L)\), is called the \emph{Weyl group} of~\(L\).
From the definition of~\(\rho_v\) it follows that \(\Weyl(L)\) is always contained in the kernel of \(\OO(L) \rightarrow \OO(L^\sharp)\).

We write \(L(n)\) for the lattice with the same underlying \(\IZ\)-module whose Gram matrix is \(nA\), where \(A\) is any Gram matrix of \(L\). 
The standard negative definite ADE lattices are denoted \(\bA_n,\bD_n,\bE_n\).

A \emph{genus} is the set of isomorphism classes of all lattices of fixed signature and discriminant form. A genus is always a finite set (see for instance \cite[Kapitel VII, Satz (21.3)]{Kneser:quadratische.Formen}).

\subsection{Preliminary results} \label{sec:preliminary_results}

Note any embedding \(\iota \colon \bU \hookrightarrow L\) is primitive, because the the lattice \((\iota(\bU) \otimes \IQ) \cap L\) is an overlattice of \(\bU\) and each overlattice of \(\bU\) is trivial (cf. \cite[Prop.~1.4.1]{Nikulin:int.sym.bilinear.forms}).

Let \(e,f\) be a fixed basis of \(\bU\) such that \(e^2 = f^2 = 0\) and \(e \cdot f = 1\).

\begin{definition} \label{def:geometric.embedding}
We say that an embedding \(\iota\colon \bU \hookrightarrow S_X\) is \emph{geometric} if \(\iota(e)\) is the class of an elliptic curve \(E\) and and \(\iota(f-e)\) is the class of a smooth rational curve \(O\) with \(E \cdot O = 1\). (Such embeddings were called ``canonical'' by Bertin et al.~ \cite{Bertin.Garbagnati.Hortsch.Lecacheux.Mase.Salgado.Whitcher}, but we believe this word to be slightly misleading.)
\end{definition}

Let \(\cE_X\) denote the set of geometric embeddings of \(\bU\) into \(S_X\) and let \(\aut(X)\) be the image of the homomorphism \(\Aut(X) \rightarrow \OO(S_X)\).

\begin{lemma} \label{lem:bijection.embeddings-fibrations}
The map
\[
\cE_X/\aut (X) \rightarrow \cJ_X/\Aut (X)
\]
defined by sending a geometric embedding \(\iota \colon \bU \hookrightarrow S_X\) to the fibration induced by the elliptic curve~\(E \coloneqq \iota(e)\) is a bijection.
\end{lemma}
\proof 
The map is clearly well defined and surjective.
Consider now two geometric embeddings \(\iota_1,\iota_2\) such that \(\iota_1(e) = \iota_2(e)\) is the class of \(E\) and suppose that \(\iota_1(f-e)\), \(\iota_2(f-e)\) are the classes of the curves \(O_1\), \(O_2\). Then, translation by a suitable section induces an automorphism \(\alpha \in \Aut(X)\) such that \(\alpha(E) = E\) and \(\alpha(O_1) = O_2\) (cf. for instance~\cite[§7.6]{Schuett.Shioda:elliptic.surfaces}). 
Therefore, \(\iota_1\) and \(\iota_2\) belong to the same \(\aut(X)\)-orbit, so the map is also injective.
\endproof

The \emph{positive cone} \(P_X\) is the connected component of
\[
    \Set{x \in S_X \otimes \IR}{x^2 > 0}
\]
that contains one ample class. 
The \emph{nef cone} \(N_X\) is defined as
\[
    N_X \coloneqq \Set{x \in S_X \otimes \IR}{x \cdot C \geq 0 \text{ for all curves } C \subset X}.
\]
Furthermore, we set
\begin{align*} 
    \OO(S_X,P_X) & \coloneqq \Set{\gamma \in \OO(S_X)}{\gamma(P_X) \subseteq P_X}, \\
    \OO(S_X,N_X) & \coloneqq \Set{\gamma \in \OO(S_X)}{\gamma(N_X) \subseteq N_X}.
\end{align*}

\begin{lemma} \label{lem:S_X}
If \(\cJ_X \neq \emptyset\), then the Néron--Severi lattice \(S_X\) is unique in its genus and the restriction of the natural homomorphism \(\OO(S_X) \rightarrow \OO(S_X^\sharp)\) to \(\OO(S_X,N_X)\) is surjective.
\end{lemma} 
\proof
If \(\cJ_X \neq \emptyset\), then there exists a (geometric) embedding \(\bU \hookrightarrow S_X\), hence \(S_X \cong \bU \oplus W\) for some lattice \(W\). 
From Nikulin's \cite[Thm.~1.14.2]{Nikulin:int.sym.bilinear.forms} we infer that \(S_X\) is unique in its genus and that the natural homomorphism \(\OO(S_X) \rightarrow \OO(S_X^\sharp)\) is surjective.

The isometry \(\gamma \in \OO(S_X)\) defined as \((-\id_\bU,\id_W)\) on the decomposition \(S_X \cong \bU \oplus W\) does not belong to \(\OO(S_X,P_X)\). Hence, \(\OO(S_X)\) is generated by \(\OO(S_X,P_X)\) and \(\gamma\).
Since \(\gamma\) is contained in the kernel of \(\OO(S_X) \rightarrow \OO(S_X^\sharp)\), the restriction of \(\OO(S_X) \rightarrow \OO(S_X^\sharp)\) to \(\OO(S_X,P_X)\) is surjective.

Moreover, it holds \(\OO(S_X,P_X) \cong \Weyl(S_X) \rtimes \OO(S_X,N_X)\) (see for instance \cite[Prop.~1.3]{Ohashi:number.Enriques}) and since \(\Weyl(S_X)\) is contained in the kernel of \(\OO(S_X) \rightarrow \OO(S_X^\sharp)\), the claim follows.
\endproof

\begin{proposition} \label{prop:effective-->section}
Let \(X \rightarrow \IP^1\) be an elliptic fibration with fiber \(E\).
If \(D\) is a divisor such that \(D \cdot E = 1\) and \(D^2 = -2\), then there exists a section \(O\) and an element of the Weyl group \(\rho \in \Weyl(S_X)\) such that 
\[
    \rho(E) = E \quad \text{and} \quad \rho(D) = O.
\]
\end{proposition}
\proof
Imitating Kond\={o}'s proof of \cite[Lemma 2.1]{Kondo:automorphisms.K3.which.act.trivially.on.Pic}, we can show that there exist a section \(O\) and \(m_0,m_1,\ldots,m_n \in \IZ\) such that
\begin{equation} \label{eq:D} 
    D = O + m_0 E + \sum_{i = 1}^n m_i C_i,
\end{equation} 
where \(C_1,\ldots,C_n\) are the irreducible fiber components such that \(O \cdot C_i = 0\) and \(C_i^2 = -2\).

From \(D^2 = -2\) it follows that
\begin{equation} \label{eq:m_0}
    m_0 = -\frac 12 \Big(\sum_{i=1}^n m_i C_i\Big)^2.
\end{equation} 
As \(C_1,\ldots,C_n\) generate a negative definite lattice, it holds \(m_0 \geq 0\). Moreover, \(m_0 = 0\) if and only if \(m_1 = \ldots = m_n = 0\). 
We claim that whenever \(m_0 > 0\) we can find \(\rho' \in \Weyl(S_X)\) such that \(\rho'(E) = E\) and 
\[
    \rho'(D) = O + m_0' E + \sum_{i = 1}^n m_i' C_i,
\]
with \(m_0' < m_0\). In a finite number of steps we obtain \(\rho = \ldots \circ \rho'' \circ \rho' \in \Weyl(S_X)\) with \(\rho(E) = E\) and \(\rho(D) = O\), thus proving the theorem. From now on, we assume that \(m_0 > 0\).

Without loss of generality, we can assume that \(C_1,\ldots,C_n\) are all components of the same fiber (otherwise we apply the following procedure on each fiber). 
Let \(C_0\) be the fiber component with \(O \cdot C_0 = 1\) and let \(\rho_i \in \Weyl(S_X)\) be the involution induced by the class of \(C_i\), for \(i = 0,1,\ldots,n\).

Note that by applying \(\rho_i\) to \(D\), with \(i \in \set{1,\ldots,n}\), only the coefficient of \(C_i\) in \eqref{eq:D} changes. In particular, the coefficients of \(O\) and \(E\) remain equal to \(1\) and \(m_0\), respectively. 
Since \(C_1,\ldots,C_n\) generate a negative definite lattice, there exist only a finite number of \(n\)-tuples \((m_1,\ldots,m_n)\) such that \eqref{eq:m_0} holds. 
If the \(n\)-tuple \((m_1,\ldots,m_{j-1},m_j',m_{j+1},\ldots,m_n)\) corresponding to \(\rho_j(D)\) satisfies \(m_j' < m_j\), we substitute \(D\) with \(\rho_j(D)\). 
Repeating this process in a finite number of steps, we can assume (up to substituting \(D\) with \(\tilde\rho(D)\) for some \(\tilde\rho \in \Weyl(S_X)\) contained in the subgroup generated by \(\rho_1,\ldots,\rho_n\)) that \(D\) satisfies the following minimality property: for each \(j = 1,\ldots,n\) it holds
\begin{equation} \label{eq:minimality_property}
    \rho_j(D) = O + m_0 E + \sum_{i = 1}^{j-1} m_i C_i + m_j' C_j + \sum_{i = j+1}^{n} m_i C_i, \quad \text{with \(m_j' \geq m_j\)}.
\end{equation}

We claim now that the coefficient \(m_0'\) of \(\rho_0(D)\) satisfies \(m_0' < m_0\), thus concluding the proof. 
We need to divide the proof according to the dual graph of \(C_1,\ldots,C_n\).

Assume first that \(n = 1\), so that the dual graph of \(C_1,\ldots,C_n\) is \(\bA_1\).
From the minimality property \eqref{eq:minimality_property} and the assumption \(m_0 > 0\) it follows that
\(m_1 < 0\).
Using \(C_0 = E - C_1\) and \(C_0 \cdot C_1 = 2\) we obtain
\[
    \rho_0(D) = (O + C_0) + m_0E + m_1(2C_0 - C_1) = O + m_0'E + m_1'C_1.
\]
with \(m_0' = m_0 + 1 + 2m_1 < m_0\), as wished.

Assume now that the dual graph of \(C_1,\ldots,C_n\) is \(\bA_n\), with \(n \geq 2\).
\[
    \dynkin[edge length=1cm,labels={1,2,n-1,n},label macro/.code={\drlap{#1}}]A{}
\]
From the minimality property \eqref{eq:minimality_property} it follows that
\[
    2m_1 \leq m_2, \qquad 
    2m_i \leq m_{i-1} + m_{i+1}, \quad \text{for } i = 2,\ldots,n-1, \qquad 
    2m_{n} \leq m_{n-1}.
\]
It holds \(im_{i-1} \leq (i-1)m_i\) for \(i = 2,\ldots,n\). Indeed, this is clear for \(i = 2\) and it follows by induction from
\[
    2im_i \leq im_{i-1} + im_{i+1} \leq (i-1)m_i + im_{i+1}.
\]
From \(nm_{n-1} \leq (n-1)m_n\) and \(2m_n \leq m_{n-1}\) we infer that \(m_n \leq 0\). 
It cannot be \(m_n = 0\) because it implies \(m_{n-1} = \ldots = m_1 = 0\), contradicting the fact that \(m_0 > 0\). 
Therefore, it holds \(m_n < 0\) and, symmetrically, \(m_1 < 0\). Using \(C_0 = E - C_1 - \ldots - C_n\), \(C_0 \cdot C_1 = C_0 \cdot C_n = 1\) and \(C_0 \cdot C_i = 0\) for \(i = 2,\ldots,n-1\) we obtain
\[
    \rho_0(D) = (O + C_0) + m_0E + m_1(C_0 - C_1) + \sum_{i = 2}^{n-1} m_iC_i + m_n(C_0 - C_n) 
        = O + m_0'E + \sum_{i = 2}^{n-1} m_i'C_i.
\]
with \(m_0' = m_0 + 1 + m_1 + m_n < m_0\), as wished.

Assume now that the dual graph of \(C_1,\ldots,C_n\) is \(\bD_n\) (\(n \geq 4\)).
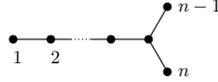
\[
    \dynkin[edge length=0.5cm,labels={1,2,{},{},n-1,n},label macro/.code={\drlap{#1}}]D{}
\]
From the minimality property \eqref{eq:minimality_property} it follows that
\begin{align*} 
    2m_1 & \leq m_2, 
    & 2m_i & \leq m_{i-1} + m_{i+1}, \quad \text{for } i = 2,\ldots,n-3, \\
    2m_{n-2} & \leq m_{n-3} + m_{n-1} + m_{n} &
    2m_{n-1} & \leq m_{n-2}, \qquad 2m_{n} \leq m_{n-2}.
\end{align*}
It holds \(m_{n-i} \leq m_{n-i-1}\) for all \(i = 2,\ldots,n-2\). Indeed, this follows from the last three inequalities for \(i = 2\) and then by induction from
\[
    2m_{n-i-1} \leq m_{n-i-2} + m_{n-i} \leq m_{n-i-2} + m_{n-i-1}.
\]
In particular \(m_2 \leq m_1\), so from \(2m_1 \leq m_2\) we infer that \(m_2 \leq 0\). 
It cannot be \(m_2 = 0\) because it implies \(m_1 = \ldots = m_n = 0\), contradicting the fact that \(m_0 > 0\). 
Moreover, it cannot be \(m_2 = -1\), because (recalling that \(m_i \in \IZ\)) it implies \(m_1 = \cdots = m_{n-2} = -1\), \(m_{n-1} \leq -1\), \(m_{n} \leq -1\), leading to the contradiction
\[
    -2 = 2m_{n-2} \leq m_{n-3} + m_{n-1} + m_n \leq -3.
\]
Therefore, it holds \(m_2 < -1\). Using \(C_0 = E - C_1 - 2C_2 - \ldots - 2C_{n-2} - C_{n-1} - C_n\), \(C_0 \cdot C_2 = 1\) and \(C_0 \cdot C_i = 0\) for \(i = 1,3,4,\ldots,n\), we obtain
\[
    \rho_0(D) = (O + C_0) + m_0E + m_1C_1 + m_2(C_0 - C_2) + \sum_{i = 3}^n m_i C_i
        = O + m_0'E + \sum_{i = 1}^n m_i'C_i.
\]
with \(m_0' = m_0 + 1 + m_2 < m_0\), as wished.

Assume now that the dual graph of \(C_1,\ldots,C_n\) is \(\bE_6\).
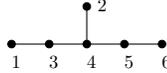
\[
    \dynkin[edge length=0.5cm,labels={1,2,3,4,5,6},label macro/.code={\drlap{#1}}]E{6}
\]
From the minimality property~\eqref{eq:minimality_property} it follows that
\begin{align*}
    2m_1 & \leq m_3, & 
    2m_2 & \leq m_4, & 
    2m_3 & \leq m_1 + m_4, \\
    2m_4 & \leq m_2 + m_3 + m_5, &
    2m_5 & \leq m_4 + m_6, &
    2m_6 & \leq m_5.
\end{align*}
From the first and third inequality we obtain \(3m_1 \leq m_4\) and from
\[
    6m_3 \leq 3m_1 + 3m_4 \leq 4m_4
\]
we get \(3m_3 \leq 2m_4\). Symmetrically, it holds \(3m_5 \leq 2m_4\). Then, from
\[
    6m_4 \leq 3m_2 + 3m_3 + 3m_5 \leq 3m_2 + 4m_4 
\]
we infer \(2m_4 \leq 3m_2\). 
Together with \(2m_2\leq m_4\), this implies \(m_2 \leq 0\). 
It cannot be \(m_2 = 0\), because it implies \(m_1 = \ldots = m_n = 0\), contradicting the fact that \(m_0 > 0\). Moreover, it cannot be \(m_2 = -1\), because (always recalling that \(m_i \in \IZ\)) it implies \(m_4 = -2\), \(m_3 \leq -2\) and \(m_5 \leq -2\), leading to the contradiction
\[
    -4 = 2m_4 \leq m_2 + m_3 + m_5 \leq -1 -2 -2 = -5.
\]
Therefore, it holds \(m_2 < -1\). Using \(C_0 = E - C_1 - 2C_2 - 2C_3 - 3C_4 - 2C_5 - C_6\), \(C_0 \cdot C_2 = 1\) and \(C_0 \cdot C_i = 0\) for \(i = 1,3,\ldots,6\), we obtain
\[
    \rho_0(D) = (O + C_0) + m_0E + m_1C_1 + m_2(C_0 - C_2) + \sum_{i = 3}^6 m_i C_i
        = O + m_0'E + \sum_{i = 1}^{6} m_i'C_i.
\]
with \(m_0' = m_0 + 1 + m_2 < m_0\), as wished.

Assume now that the dual graph of \(C_1,\ldots,C_n\) is \(\bE_7\).
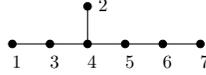
\[
    \dynkin[edge length=0.5cm,labels={1,2,3,4,5,6,7},label macro/.code={\drlap{#1}}]E{7}
\]
From the minimality property~\eqref{eq:minimality_property} it follows that
\begin{align*}
    2m_1 & \leq m_3, & 
    2m_2 & \leq m_4, & 
    2m_3 & \leq m_1 + m_4, &
    2m_4 & \leq m_2 + m_3 + m_5, \\
    2m_5 & \leq m_4 + m_6, &
    2m_6 & \leq m_5 + m_7, &
    2m_7 & \leq m_6.
\end{align*}
From the last three inequalities we obtain \(3m_6 \leq 2m_5\) and \(4m_5 \leq 3m_4\). From
\[
    8m_4 \leq 4m_2 + 4m_3 + 4m_5 \leq 2m_4 + 4m_3 + 3m_4
\]
we get \(3m_4 \leq 4m_3\). Then, from
\[
    6m_3 \leq 3m_1 + 3m_4 \leq 3m_1 + 4m_3  
\]
we infer \(2m_3 \leq 3m_1\). Together with \(2m_1 \leq m_3\), this implies \(m_1 \leq 0\).
It cannot be \(m_1 = 0\), because it implies \(m_2 = \ldots = m_7 = 0\), contradicting the fact that \(m_0 > 0\). Moreover, it cannot be \(m_1 = -1\), because (recalling that \(m_i \in \IZ\)) it implies \(m_3 = -2\), \(m_4 = -3\), \(m_2 \leq -2\), \(m_5 \leq -3\), and \(m_5 \leq -2\), leading to the contradiction
\[
    -6 = 2m_4 \leq m_2 + m_3 + m_5 \leq -2 -2 -3 = -7.
\]

Therefore, it holds \(m_1 < -1\). Using \(C_0 = E - 2C_1 - 2C_2 - 3C_3 - 4C_4 - 3C_5 - 2C_6-C_7\), \(C_0 \cdot C_1 = 1\) and \(C_0 \cdot C_i = 0\) for \(i = 2,\ldots,7\), we obtain
\[
    \rho_0(D) = (O + C_0) + m_0E + m_1(C_0 - C_1) + \sum_{i = 1}^7 m_i C_i
        = O + m_0'E + \sum_{i = 2}^{6} m_i'C_i.
\]
with \(m_0' = m_0 + 1 + m_1 < m_0\), as wished.

Finally, assume that the dual graph of \(C_1,\ldots,C_n\) is \(\bE_8\).
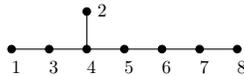
\[
    \dynkin[edge length=0.5cm,labels={1,2,3,4,5,6,7,8},label macro/.code={\drlap{#1}}]E{8}
\]
From the minimality property~\eqref{eq:minimality_property} it follows that
\begin{align*}
    2m_1 & \leq m_3, & 
    2m_2 & \leq m_4, & 
    2m_3 & \leq m_1 + m_4, &
    2m_4 & \leq m_2 + m_3 + m_5, \\
    2m_5 & \leq m_4 + m_6, &
    2m_6 & \leq m_5 + m_7, &
    2m_7 & \leq m_6 + m_8, &
    2m_8 & \leq m_7.
\end{align*}
As in the \(\bE_6\) case, we obtain \(3m_3 \leq 2m_4\). From
\[
   12m_4 \leq 6m_2 + 6m_3 + 6m_5 \leq 3m_4 + 4m_3 + 6 m_5  
\]
we get \(5m_4 \leq 6m_5\). Then, from
\[
    10m_5 \leq 5m_4 + 5m_6 \leq 6m_5 + 5m_6
\]
we get \(4m_5 \leq 5m_6\). 
Similarly, we infer \(3m_6 \leq 4m_7\) and \(2m_7 \leq 3m_8\).
Together with \(2m_8 \leq m_7\), this implies \(m_8 \leq 0\). 
It cannot be \(m_8 = 0\), because it implies \(m_7 = \ldots = m_1 = 0\), contradicting the fact that \(m_0 > 0\). 
Moreover, it cannot be \(m_8 = -1\), because (recalling that \(m_i \in \IZ\)) it implies \(m_7 = -2\), \(m_6 = -3\), \(m_5 = -4\), \(m_4 = -5\), \(m_3 \leq -4\) and \(m_2 \leq -3\), leading to the contradiction
\[
    -10 = 2m_4 \leq m_2 + m_3 + m_5 \leq -3 - 4 - 4 = -11.
\]
Therefore, it holds \(m_8 < -1\). Using \(C_0 = E - 2C_1 - 3C_2 - 4C_3 - 6C_4 - 5C_5 - 4C_6 - 3C_7 - 2C_8\), \(C_0 \cdot C_8 = 1\) and \(C_0 \cdot C_i = 0\) for \(i = 1,\ldots,7\), we obtain
\[
    \rho_0(D) = (O + C_0) + m_0E + \sum_{i = 1}^7 m_i C_i + m_8(C_0 - C_8) 
        = O + m_0'E + \sum_{i = 1}^{8} m_i'C_i.
\]
with \(m_0' = m_0 + 1 + m_8 < m_0\), as wished.
\endproof

\begin{definition} \label{def:frame_genus}
The \emph{frame genus} of \(X\) is the genus~\(\cW_X\) of negative definite lattices \(W\) such that
\[
    \rank(W) = \rank(S_X) - 2 \quad \text{and} \quad W^\sharp \cong S_X^\sharp.
\]
We define the \emph{frame map} \(\fr_X\) as follows:
\[ 
    \fr_X \colon \cE_X/{\aut(X)} \rightarrow \cW_X, \qquad \iota \mapsto \iota(\bU)^\perp.
\]
The number \(|{\fr_X^{-1}(W)}|\) is called the \emph{multiplicity} of the frame~\(W \in \cW_X\).
\end{definition}

\begin{corollary} \label{lem:surjective}
For each frame \(W \in \cW_X\) there exists a geometric embedding \(\iota \colon \bU \hookrightarrow S_X\) such that \(\iota(\bU)^\perp \cong W\). In particular, the frame map \(\fr_X\) is surjective and it holds
\begin{equation}  \label{eq:|orbits|=sum(multiplicities)}
    |\cJ_X/{\Aut(X)}| = \sum_{W \in \cW_X} |{\fr_X^{-1}(W)}|.
\end{equation}
\end{corollary} 
\proof
Fix a frame \(W \in \cW_X\). The lattices \(S_X\) and \(\bU \oplus W\) belong to the same genus. Since \(S_X\) is unique in its genus by \autoref{lem:S_X}, it holds \(S_X \cong \bU \oplus W\), i.e. there exists an embedding \(\iota \colon \bU \hookrightarrow S_X\) with \(\iota(\bU)^\perp \cong W\). 

Since \(\Weyl(S_X)\) acts transitively on the chambers of the positive cone, we find \(\rho' \in \Weyl(S_X)\) such that \(\rho'\circ \iota(e) \in N_X\), hence \(\rho'\circ \iota(e)\) is the class of an elliptic curve~\(E\) (see \cite[§3, proof of Cor.~3]{Pjateckii-Shapiro.Shafarevich:Torelli's.theorem.K3} or \cite[Ch.~2, Prop.~3.10]{Huybrechts:lectures.K3}).
Then, \(D = \rho'\circ \iota(f-e)\) satisfies the hypothesis of~\autoref{prop:effective-->section}, so there exists \(\rho \in \Weyl(S_X)\) such that \(\rho \circ \rho' \circ \iota\) is a geometric embedding.
Clearly, it holds \((\rho\circ\rho' \circ \iota(\bU))^\perp \cong \iota(\bU)^\perp \cong W\).

By \autoref{lem:bijection.embeddings-fibrations} we have \(|\cJ_X/{\Aut(X)}| = |\cE_X/\aut(X)|\), from which we infer equation~\eqref{eq:|orbits|=sum(multiplicities)}.
\endproof 

\begin{proposition} \label{prop:ample}
Let \(\pi\colon X \rightarrow \IP^1\) be an elliptic fibration with fiber \(E\) and section \(O\).
If \(C_1,\ldots,C_n\) are the components of the reducible fibers not intersecting~\(O\), then there exist \(n_0, m_0,\ldots,m_n \in \IZ\) such that 
\(
    D = n_0 O + m_0 E + \sum_{i = 1}^n m_i C_i
\)
is an ample divisor.
\end{proposition}
\proof
For the sake of simplicity, we assume that there is only one reducible fiber, but the same argument works if there is more than one.
Let \(C_0\) be the rational fiber component with \(O \cdot C_0 = 1\).
The sublattice generated by \(C_1,\ldots,C_n\) is of ADE type, so we can find \(R \coloneqq \sum_{i = 1}^n m_i C_i\) in the suitable Weyl chamber  in such a way that
\[
    R \cdot C_i > 0 \qquad \text{for each \(i = 1,\ldots n\)}.
\]

Fix a positive \(n_0 \in \IZ\) such that \(n_0 > -R \cdot C_0\). Since \(E\) is linearly equivalent to a certain combination of \(C_0,\ldots,C_n\) with positive coefficients, taking \(m_0\) large enough we can suppose that \(D = n_0 O + m_0 E + R\) is linearly equivalent to
\[
    n_0 O + m_0 C_0 + \sum m_i' C_i \qquad \text{with \(m_i' > 0\)}.
\]

Let \(C\) be a smooth rational curve on \(X\). If \(C \neq O\) and \(C \neq C_i\) for each \(i = 0,\ldots,n\), then \(D \cdot C > 0\) because \(C\) must intersect one of the components \(C_i\). 

If \(C = C_i\) for some \(i \in \set{1,\ldots,n}\), then
\[
    D \cdot C = R \cdot C_i > 0. 
\]
If \(C = C_0\), then
\[
    D \cdot C = n_0 + R \cdot C_0 > 0.
\]
If \(C = O\), then
\[
    D \cdot C = -2n_0 + m_0,
\]
which is positive if \(m_0\) is large enough. Finally, it holds
\[
    D^2 = -2n_0^2 + 2n_0m_0 + R^2.
\]
Therefore, if \(m_0\) is large enough, then \(D^2 > 0\) and \(D \cdot C > 0\) for every smooth rational curve on~\(X\), proving that \(D\) is ample (cf. for instance \cite[Ch.~8, Cor.~1.6]{Huybrechts:lectures.K3}).
\endproof

\subsection{Statement and corollaries} \label{sec:statement}

The Hodge decomposition on \(\HH^2(X,\IZ)\) induces a Hodge decomposition on the transcendental lattice~\(T_X\). 
The group of Hodge isometries of~\(T_X\) is denoted \(\OO_\hdg(T_X)\). 
If \(\sigma_X\) is a generator of the subspace \(\HH^{2,0}(X) \subset T_X \otimes \IC\), we have
\[
    \OO_\hdg(T_X) \coloneqq \Set{\eta \in \OO(T_X)}{\eta(\sigma_X) \in \IC \sigma_X}.
\]
The image of the natural homomorphism
\[
\OO_\hdg(T_X) \rightarrow \OO(T_X^\sharp)
\]
is denoted \(\OO_\hdg^\sharp(T_X)\). 
We are now able to state the main theorem, whose proof is contained in~\autoref{sec:thm:main_mutliplicity}.

\begin{theorem} \label{thm:main_multiplicity}
Let \(X\) be a complex projective K3 surface with frame genus~\(\cW_X\) (\autoref{def:frame_genus}). Then, the multiplicity of a frame \(W \in \cW_X\) is given by
\[
    |{\fr_X^{-1}(W)}| = |{\OO_\hodge^\sharp(T_X)} \backslash {\OO(T_X^\sharp)} / {\OO^\sharp(W)}|.
\]
\end{theorem}

\begin{remark} \label{rem:Cauchy-Frobenius}
By~\autoref{def:frame_genus} it holds \(W^\sharp \cong S_X^\sharp \cong T_X(-1)^\sharp\) for each \(W \in \cW_X\), so the group~\(\OO^\sharp(W)\) can be considered as a subgroup of \(\OO(T_X^\sharp)\) thanks to the (non-canonical) isomorphisms
\[
    \OO(W^\sharp) \cong \OO(S_X^\sharp) \cong \OO(T_X^\sharp),
\]
which we fix once and for all.
The subgroup \(\OO^\sharp(W)\), therefore, is only well-defined up to conjugation. 
Still, the formula of \autoref{thm:main_multiplicity} makes sense because the number of double cosets \(|H \backslash G / K|\) for any subgroups \(H,K\) of a group \(G\) only depends on the conjugacy class of \(H\) and \(K\). 
Indeed, we have the following formula due to Cauchy and Frobenius:
\begin{equation} \label{eq:Cauchy-Frobenius}
    |H \backslash G / K| = \frac{1}{|H||K|} \sum_{(h,k) \in H \times K} |{\Set{g \in G}{hgk = g}}|.
\end{equation}
\end{remark}


The following corollary providing uniform bounds on the number of orbits \(|\cJ_X/{\Aut(X)}|\) already appears in previous literature under various forms, cf. for instance {\cite[Prop.~C']{Braun.Kimura.Watari:classif.ell.fibr}} and {\cite[Thm~3.10]{Mezzedimi:K3.zero.entropy}}.

\begin{corollary}  \label{cor:uniform.bound}
It always holds
\[
    |\cW_X| \leq |\cJ_X/{\Aut(X)}| \leq |\cW_X|\cdot |{\OO_\hdg^\sharp(T_X)}\backslash{\OO(T_X^\sharp)}|.
\]
In particular, \(|\cJ_X/{\Aut(X)}| = |\cW_X|\) when the map \(\OO_\hdg(T_X)\rightarrow \OO(T_X^\sharp)\) is surjective.
\end{corollary}
\proof
The map \(\fr_X\) is surjective by \autoref{lem:surjective}, hence \(|{\fr_X^{-1}(W)}| \geq 1\) for each \(W \in \cW_X\). By~\eqref{eq:|orbits|=sum(multiplicities)},
\[
    |\cW_X| \leq \sum_{W \in \cW_X} |{\fr_X^{-1}(W)}| = |\cJ_X/{\Aut(X)}|,
\]
proving the first inequality. Moreover, it holds trivially
\[
    |{\OO_\hdg^\sharp(T_X)}\backslash{\OO(T_X^\sharp)}/\OO^\sharp(W)| \leq |{\OO_\hdg^\sharp(T_X)}\backslash{\OO(T_X^\sharp)}|.
\]
By \autoref{thm:main_multiplicity} and again equation~\eqref{eq:|orbits|=sum(multiplicities)} we obtain therefore
\[
    |\cJ_X/{\Aut(X)}| = \sum_{W \in \cW_X} |{\fr_X^{-1}(W)}| \leq \sum_{W \in \cW_X} |{\OO_\hdg^\sharp(T_X)}\backslash{\OO(T_X^\sharp)}| = |\cW_X| \cdot |{\OO_\hdg^\sharp(T_X)}\backslash{\OO(T_X^\sharp)}|. \qedhere 
\]
\endproof

\subsection{Proof of \autoref{thm:main_multiplicity}} \label{sec:thm:main_mutliplicity}

Let \(G \coloneqq \OO(T_X^\sharp)\) and \(H \coloneqq \OO_\hdg^\sharp(T_X)\).
Fix a frame \(W \in \cW_X\) and a geometric embedding \(\iota_0 \in \fr_X^{-1}(W)\) (\autoref{def:geometric.embedding}), which exists by \autoref{lem:surjective}. Since \(\bU\) is unimodular, we have a decomposition
\[
    S_X \cong \iota_0(\bU) \oplus \iota_0(\bU)^\perp \cong \bU \oplus W.
\]

Consider the group \(\OO(S_X, \iota_0)\) of isometries \(S_X\) preserving this decomposition:
\[
    \OO(S_X, \iota_0) \coloneqq \Set{\gamma \in \OO(S_X)}{\gamma \circ \iota_0(\bU) = \iota_0(\bU)}.
\]
Let \(K\) be the image of \(\OO(S_X, \iota_0)\) in \(G\). Using the isomorphism \(S_X^\sharp \cong \bU^\sharp \oplus W^\sharp \cong W^\sharp\) induced by this decomposition, the subgroup~\(K\) can be identified with \(\OO^\sharp(W)\) (cf. \autoref{rem:Cauchy-Frobenius}).

Recall that the Torelli theorem for K3 surfaces~\cite{Pjateckii-Shapiro.Shafarevich:Torelli's.theorem.K3} asserts that
\begin{equation} \label{eq:Torelli}
    \aut(X) = \Set{ \gamma \in \OO(S_X)}{\gamma \in \OO(S_X,N_X) \text{ and } \gamma^\sharp \in H}.
\end{equation}

We define a map
\[
    \dc \colon \fr_X^{-1}(W) \rightarrow H \backslash G / K
\]
in the following way. Take a geometric embedding \(\iota \in \fr_X^{-1}(W)\). 
Recalling that \(S_X\) is unique in its genus by \autoref{lem:S_X}, we see that the embedding is unique up to isometries of~\(S_X\) thanks to Nikulin's \cite[Prop.~1.15.1]{Nikulin:int.sym.bilinear.forms}. 
This means that there exists \(\gamma \in \OO(S_X)\) such that \(\iota = \gamma \circ \iota_0\). 
Using the identification \(\OO(T_X^\sharp) \cong \OO(S_X^\sharp)\), we set
\[
    \dc(\iota) \coloneqq H \gamma^\sharp K.
\]

We claim that the map \(\dc\) is well defined. 
Indeed, take two geometric embeddings \(\iota_1,\iota_2 \in \fr_X^{-1}(W)\) such that \(\iota_2 = \alpha \circ \iota_1\) for some \(\alpha \in \aut(X)\).  
Let \(\iota_1 = \gamma_1 \circ \iota_0\) and \(\iota_2 = \gamma_2 \circ \iota_0\) for some \(\gamma_1,\gamma_2 \in \OO(S_X)\).
The isometry \(\kappa \coloneqq \gamma_1^{-1} \circ \alpha^{-1} \circ \gamma_2\) belongs to \(\OO(S_X, \iota_0)\) because
\[
    \kappa \circ \iota_0 = \gamma_1^{-1} \circ \alpha^{-1} \circ (\gamma_2 \circ \iota_0) = \gamma_1^{-1} \circ (\alpha^{-1} \circ \iota_2) = \gamma_1^{-1} \circ \iota_1 = \iota_0.
\]
In particular, it holds \(\kappa^\sharp \in K\). 
Since \(\gamma_2 = \alpha \circ \gamma_1 \circ \kappa\) and \(\alpha^\sharp \in H\)
 by \eqref{eq:Torelli}, the elements \(\gamma_1^\sharp,\gamma_2^\sharp \in G\) belong to the same \((H,K)\)-double coset.

Next we claim that the map \(\dc\) is injective. 
Indeed, take two geometric embeddings \(\iota_1 = \gamma_1\circ \iota_0\), \(\iota_2 = \gamma_2 \circ \iota_0\), with \(\gamma_1,\gamma_2 \in \OO(S_X)\), such that \(\gamma_1^\sharp,\gamma_2^\sharp\) belong to the same \((H,K)\)-double coset.
By definition of~\(H\) and \(K\), this means that there exist \(\eta \in \OO_\hdg(T_X)\) and \(\kappa \in \OO(S_X, \iota_0)\) such that \(\gamma_2^\sharp = \eta^\sharp \gamma_1^\sharp \kappa^\sharp\).
We need to show that there exists \(\alpha \in \aut(X)\) such that \(\iota_2 = \alpha \circ \iota_1\).

We define \(\alpha_U \coloneqq \gamma_2 \circ \gamma_1^{-1}\) on \(\iota_1(\bU)\) and \(\alpha_W \coloneqq \gamma_2 \circ \kappa^{-1} \circ \gamma_1^{-1}\) on \(\iota_1(\bU)^\perp\). 
As \(\bU\) is unimodular, \(\alpha_U\) and \(\alpha_W\) glue well together on \(\iota_1(\bU) \oplus \iota_1(\bU)^\perp \cong S_X\), defining an element \(\alpha\in \OO (S_X)\).
From the definition, it immediately follows that 
\[
    \iota_2 = \gamma_2 \circ \iota_0 = \alpha \circ \gamma_1 \circ \iota_0 = \alpha \circ \iota_1 \qquad \text{and} \qquad \alpha^\sharp = (\gamma_2 \circ \kappa^{-1} \circ \gamma_1^{-1})^\sharp = \eta^\sharp \in H.
\]
For \(i = 1,2\) consider the fibrations \(\pi_i\colon X \rightarrow \IP^1\) corresponding to \(\iota_i\). Write \(E_i \coloneqq \iota_i(e)\), \(O_i \coloneqq \iota_i(f-e)\) and let \(C_{i1},\ldots,C_{in}\) be the irreducible components of the reducible fibers of \(\pi_i\). 
Thanks to \autoref{prop:ample} we know that there exist \(n_0,m_0,\ldots,m_n\) such that \(D = n_0 O_1 + m_0 E_1 + \sum_{i=1}^n m_i C_{1i}\) is an ample divisor. Set 
\[
    R \coloneqq \alpha(D) - n_0 O_2 - m_0 E_2 = \sum_{i=1}^n m_i \alpha(C_{1i}).
\]
As each \(\alpha(C_{1i})\) is a vector of square \(-2\) in \(\iota_2(\bU)^\perp\), there exist \(m_1',\ldots,m_n' \in \IZ\) such that
\[
    R = \sum_{i=1}^nm_i' C_{2i}.
\]
Since the Weyl group of the ADE lattices acts transitively on the Weyl chambers (see for instance~\cite{Hall:Lie.groups.algebras.representations}), we can find \(\rho \in \Weyl((\iota_2(\bU)^\perp)_\rootlattice)\) such that \(\rho(R) \cdot C_{2i} > 0\) for every \(i = 1,\ldots,n\). 
As \(\rho\) acts trivially on the discriminant group of \((\iota_2(\bU)^\perp)_\rootlattice\), we can extend it to an isometry \(\tilde\rho \in \OO(\iota_2(\bU)^\perp)\) with \(\tilde\rho^\sharp = 1\). 
Therefore, up to substituting \(\alpha\) with \(\tilde\rho \circ \alpha\), we can suppose that \(\alpha(D)\) is an ample divisor by the same arguments as in \autoref{prop:ample}. 
In this way, \(\alpha\) sends an ample divisor to an ample divisor, so \(\alpha \in \OO(S_X, N_X)\) because the ample cone is the interior of \(N_X\) (see for instance~\cite[Ch.~8, Cor.~1.4]{Huybrechts:lectures.K3}). 
The Torelli Theorem~\eqref{eq:Torelli} implies that \(\alpha \in \aut(X)\), as wished.

Finally, we claim that the map \(\dc\) is surjective.
Indeed, for every double coset \(HgK, g \in G\), we can assume that \(g = \gamma^\sharp\) for some \(\gamma \in \OO(S_X,N_X)\) thanks to \autoref{lem:S_X}.
As \(\iota_0(e) \in N_X\), it follows that \(\gamma \circ \iota_0(e)\) is a primitive   nef class of square \(0\), hence \(\gamma \circ \iota_0(e)\) is the class of some elliptic curve \(E\). 
The class \(D = \gamma \circ \iota_0(f-e)\) then satisfies the assumptions of \autoref{prop:effective-->section}. 
Therefore, there exists \(\rho \in \Weyl(S_X)\) and a section~\(O\) such that \(\rho(D)= O\), which means that the embedding~\(\iota = \rho \circ \gamma \circ \iota_0\) is geometric and satisfies
\[
    \dc(\iota) = H(\rho \circ \gamma)^\sharp K = H\gamma^\sharp K = HgK.
\]

Therefore, the map \(\dc\) is a bijection and the theorem follows.
\qed

\subsection{Guidelines} \label{sec:guidelines}

We now would like to explain how to compute \(|\cJ_X/{\Aut(X)}|\) in explicit cases.
According to \autoref{thm:main_multiplicity}, the first ingredient to be computed is the image \(\OO_\hdg^\sharp(T_X)\) of the natural homomorphism
\(
    \OO_\hdg(T_X) \rightarrow \OO(T_X^\sharp).
\)
The following proposition summarizes the main properties of the group \(\OO_\hdg(T_X)\).

\begin{proposition}[see {\cite[Thm.~3.1]{Nikulin:finite.aut.groups.K3} or \cite[Prop.~B.1]{Hosono.Lian.Oguiso.Yau:Fourier-Mukai.number.K3.surf}}] \label{prop:HLOY.Appendix.B}
If \(T_X\) is the transcendental lattice of a projective K3 surface \(X\), then the group \(\OO_\hdg(T_X)\) is a finite cyclic group of even order containing~\(-\id\) such that 
\[ \varphi(|{\OO_\hdg(T_X)}|) \mid \rank(T_X), \] 
where \(\varphi\) denotes Euler's totient function. The action of \(\OO_\hdg(T_X)\) on \(\IC\sigma_X\) is faithful. \qed
\end{proposition}

The last assertion means that if \(\gamma \in \OO_\hdg(T_X)\) has order \(m\), then \(\gamma(\sigma_X) = \zeta_m \sigma_X\), where \(\zeta_m\) is a primitive \(m\)th root of unity.

\begin{remark} \label{rmk:rk.TX.odd}
If \(\rank(T_X)\) is odd, \autoref{prop:HLOY.Appendix.B} implies that \(\OO_\hdg(T_X) = \set{\pm \id}\). Therefore, \(\OO_\hdg^\sharp(T_X) = \set{\pm \id}\). 
The subgroup \(\OO^\sharp(W)\) always contains \(\set{\pm \id}\), so by \autoref{thm:main_multiplicity} we have
\[
    |{\fr_X^{-1}(W)}| = |{\OO(T_X^\sharp)}/{\OO^\sharp(W)}| = |{\OO(T_X^\sharp)}|/|{\OO^\sharp(W)}|
\]
for every \(W \in \cW_X\).
\end{remark}

A direct consequence of \autoref{prop:HLOY.Appendix.B} is that the image subgroup \(\OO_\hdg^\sharp(T_X)\) is also a finite cyclic subgroup. In all explicit examples we will be interested in determining the conjugacy class of one of its generators. 
Clearly, the order of a generator depends on the size of the kernel of the map \(\OO_\hdg(T_X) \rightarrow \OO(T_X^\sharp)\). The following lemma was already stated without proof in \cite[Rmk.~2.13]{Shimada.Veniani:Enriques.involutions.singular.K3}.

\begin{lemma} \label{lem:ker.O_hdg}
The kernel of the map \(\OO_\hdg(T_X) \rightarrow \OO(T_X^\sharp)\) is isomorphic to the kernel of the map \(\Aut(X) \rightarrow \OO(S_X)\).
\end{lemma} 
\proof 
Indeed, every element in the kernel of \(\Aut(X) \rightarrow \OO(S_X)\) induces a Hodge isometry of \(T_X\) acting trivially on \(S_X^\sharp \cong T_X^\sharp\) and,
conversely, every element in the kernel of \(\OO_\hdg(T_X) \rightarrow \OO(T_X^\sharp)\) can be extended to a Hodge isometry \(\gamma\) of \(\HH^2(X,\IZ)\) by the identity on \(S_X\). As \(\gamma\) obviously acts trivially on the nef cone, \(\gamma\) is the pullback of an automorphism of \(X\) by the Torelli Thm.~\eqref{eq:Torelli}.
\endproof 

\autoref{lem:ker.O_hdg} will prove very useful, as the classification of the possible kernels of the map \(\Aut(X) \rightarrow \OO(S_X)\) was started by Vorontsov~\cite{Vorontsov}, followed by Kond\={o}~\cite{Kondo:automorphisms.K3.which.act.trivially.on.Pic}, Oguiso and Zhang~\cite{Oguiso.Zhang:on.Vorontsov's.thm}, Schütt~\cite{Schuett:K3.non-sp.aut.2-power} and finally completed by Taki~\cite{Taki:classification.non-sp.aut.K3.which.act.trivially}.

The second and last ingredient to be computed is a list of the groups \(\OO^\sharp(W)\) for \(W \in \cW_X\). 
Therefore, one first needs to determine the genus \(\cW_X\).

Nishiyama introduced a method to compute \(\cW_X\) using a negative definite lattice \(T_0\) such that
\[
    \rank(T_0) = \rank(T_X) + 4 \quad \text{and} \quad  T_0^\sharp \cong T_X^\sharp.
\]
The method is usually known as the ``Kneser--Nishiyama method'' and consists in finding all primitive embeddings of \(T_0\) up to isometries into the 24 Niemeier lattices~\(N\). 
The orthogonal complements \(W = T_0^\perp \subset N\) run over all \(W \in \cW_X\) by \cite[Prop.~1.15.1]{Nikulin:int.sym.bilinear.forms}.
We refer to the original paper by Nishiyama~\cite{Nishiyama:Jacobian.fibrations.K3.MW} and to the surveys by Schütt and Shioda~\cite{Schuett.Shioda:elliptic.surfaces}, \cite[§12.3]{Schuett.Shioda:Mordell-Weil} for further details.

Still, since we are interested in obtaining an explicit Gram matrix of each frame \(W \in \cW_X\), we point out the following algorithm. Let \(N\) be a Niemeier lattice such that there exists a primitive embedding \(T_0 \hookrightarrow N_\rootlattice\).
\begin{itemize}
    \item Let the primitive embedding \(T_0 \hookrightarrow N_{\rootlattice}\) be given by a matrix~\(D\).
    \item Find the isotropic subgroup \(H \subset N_\rootlattice^\sharp\) corresponding to \(N_{\rootlattice} \hookrightarrow N\) (cf. \cite[Prop.~1.4.1]{Nikulin:int.sym.bilinear.forms}).
    \item Starting from \(H\), find a matrix \(C\) with rational coefficients such that \(N = CN_{\rootlattice}C^\intercal\).
    \item The embedding \(T_0 \hookrightarrow N\) is then given by the matrix \(DC^{-1}\).
    \item In this way we can compute the Gram matrix \(W\) of the orthogonal complement of \(T_0\) in \(N\), using for instance the command \verb+orthogonal_complement+ of the class \verb+IntegralLattice+ in Sage.
\end{itemize}

Once a Gram matrix is given, efficient methods to find a set of generators of \(\OO(W)\) are well known (see for instance the work by Plesken and Souvignier~\cite{Plesken.Souvignier:computing.isometries}). It is then elementary to compute the image group \(\OO^\sharp(W)\).

\begin{remark} \label{rem:Kneser.neighbor}
Another efficient method to compute a genus is Kneser's neighbor method. 
It has been implemented in several computer algebra systems (for instance, in the Magma function \verb+GenusRepresentatives+), see the work by Scharlau and Hemkemeier~\cite{Scharlau.Hemkemeier}.
\end{remark}

\begin{remark} \label{rem:mass}
The order \(|{\OO(W)}|\) and the mass of \(\cW_X\) (for the definition of \emph{mass} see the work by Conway and Sloane~\cite{Conway.Sloane:mass.formula}) can be efficiently computed using the commands \verb+number_of_automorphisms+ and \verb+conway_mass+ of the \verb+QuadraticForm+ class in Sage. 
Once a list of lattices \(W \in \cW_X\) is obtained, one can verify a posteriori that the Smith--Minkowski--Siegel mass formula holds in order to check that the list is complete.
\end{remark}

\section{Examples} \label{sec:examples}
In this section we show how to apply \autoref{thm:main_multiplicity} in order to compute \(|\cJ_X/{\Aut(X)}|\) in various explicit examples.
We first consider K3 surfaces with transcendental lattice \(T_X \cong \bU(2)^2\) in \autoref{sec:Barth-Peters}, then with \(T_X \cong \bU(2)^2\) in \autoref{sec:Oguiso}, with \(T_X \cong \bU(2)^2 \oplus [-4]\) in \autoref{sec:Kumar}, with \(T_X \cong \bU(2)^2 \oplus [-2]^2\) in \autoref{sec:Kloosterman} and finally with \(T_X \cong \bU \oplus [12]\) in \autoref{sec:Apery-Fermi}.

\subsection{Barth--Peters family} \label{sec:Barth-Peters}

Let \(X\) be a K3 surface with transcendental lattice
\[
    T_X \cong \bU \oplus \bU(2).
\]

\begin{remark}
The surface \(X\) is the generic element of a \(2\)-dimensional family introduced by Barth and Peters~\cite{Barth.Peters:aut.enriques} in order to construct Enriques surfaces with an involution acting trivially on cohomology. 
Elliptic fibrations on this family were also studied by Hulek and Schütt~\cite{Hulek.Schuett:Enriques.surf.jacobian.ell}.
\end{remark}

\begin{lemma} \label{lem:O.hdg.Barth_Peters}
It holds \(\OO_\hodge^\sharp(T_X) = \set{\id}\). 
\end{lemma}
\proof 
We choose a basis \(t_1,\ldots,t_4 \in T_X\) so that the corresponding Gram matrix is
\[
    T = {\small \begin{pmatrix}
    0 & 1 & 0 & 0 \\
    1 & 0 & 0 & 0 \\
    0 & 0 & 0 & 2 \\
    0 & 0 & 2 & 0 \\
    \end{pmatrix}}.
\]
The discriminant group~\(T_X^\sharp\) is generated by the classes of \(\frac{1}{2}t_3\) and \(\frac{1}{2}t_4\). 
The involution exchanging these two classes generates the group \(\OO(T_X^\sharp)\).

By {\cite[Theorem 3.1]{Nikulin:finite.aut.groups.K3} or \cite[Proposition~B.1]{Hosono.Lian.Oguiso.Yau:Fourier-Mukai.number.K3.surf}, it holds \(\varphi(|{\OO_\hdg(T_X)}|) \mid \rank(T_X) = 4\), therefore
\(
    |{\OO_\hdg(T_X)}| \in \set{2,4,6,8,10,12}.
\)
Since \(|{\OO(T_X^\sharp)}| = 2\), we only need to check that any isometry \(\eta \in \OO_\hodge(T_X)\) of order \(2^k\), \(k > 0\), acts trivially on \(T_X^\sharp\). 

If \(\eta\) has order~\(2\), then \(\eta = -\id\) by \cite[Proposition~B.1]{Hosono.Lian.Oguiso.Yau:Fourier-Mukai.number.K3.surf}, so \(\eta^\sharp = \id\). 
Moreover, by \autoref{lem:ker.O_hdg} and Taki's~\cite[Main theorem on page 18]{Taki:classification.non-sp.aut.K3.which.act.trivially}, \(\eta\) cannot have order \(8\). 

Hence, suppose that \(\eta\) has order \(4\). We can assume that \(\eta\) is represented by a matrix~\(A\) with \(ATA^\intercal = T\) and \(A^2 = -I\), such that the generator \(\sigma = \sigma_X\) of the subspace \(\HH^{2,0}(X) \subset T_X \otimes \IC\) is contained in the \(2\)-dimensional, totally isotropic eigenspace \(V \subset T_X \otimes \IC\) associated with \(\ri = \sqrt{-1}\). Let \(\tau \in V\) be linearly independent from \(\sigma\) and let us write
\begin{align*} 
    \sigma &= \sigma_1 t_1 + \sigma_2 t_2 + \sigma_3t_3 + \sigma_4t_4, \quad \sigma_i \in \IC, \\
    \tau &= \tau_1 t_1 + \tau_2 t_2 + \tau_3 t_3 + \tau_4 t_4 \in V, \quad \tau_i \in \IC. 
\end{align*}
It holds \(\sigma_1 \neq 0\), otherwise \(t_2 \in T_X^\perp\), so we can rescale \(\sigma_1 = 1\).
Up to substituting \(\tau\) with a linear combination of \(\sigma\) and \(\tau\) we can suppose \(\tau_1 = 0\). 
In addition, using the relations \(\sigma^2 = \tau^2 = \sigma\cdot \tau = 0\), up to substituting \(t_3\) with \(t_4\) we can assume that
\[
    \sigma_2 = -2\sigma_3 \sigma_4, \qquad \tau_2 = -2\sigma_3\tau_4, \qquad \tau_3 = 0, \qquad \tau_4 = 1.
\]

Note that \(\Im(\sigma_3) \neq 0\) because \(\sigma \cdot \bar \sigma > 0\).
Imposing that \(\sigma,\tau \in V\) and \(\bar\sigma,\bar\tau \in \bar V\), we infer that there exist \(m_1,m_2 \in \IZ\) such that
\[
    A = {\small \begin{pmatrix}
        m_1 & 0 & (m_1^2+1)/(2m_2) & 0 \\
        0 & -m_1 & 0 & m_2 \\
        -2m_2 & 0 & -m_1 & 0 \\
        0 & -(m_1^2 + 1)/m_2 & 0 & m_1 \\
    \end{pmatrix}}, \quad m_1 = -\frac{\Re(\sigma_3)}{\Im(\sigma_3)}, \quad m_2 = -\frac{1}{2\Im(\sigma_3)}.
\]
Therefore, \(A\) acts trivially on the classes of \(\frac 12 t_3\) and \(\frac 12 t_4\), so \(\eta^\sharp = \id\).
\endproof 

\begin{proposition}
The frame genus \(\cW_X\) contains exactly \(6\) isomorphism classes, listed in \autoref{tab:genus_Barth_Peters}, whose Gram matrices are contained in the arXiv ancillary file \verb+genus_Barth_Peters.sage+.
\end{proposition}
\proof 
We computed the frame genus \(\cW_X\) using the command \verb+GenusRepresentatives+ in Magma (\autoref{rem:Kneser.neighbor}). 
Alternatively, one could have applied the Kneser--Nishiyama method with \(T_0 = \bD_8\). 
The list is complete because the mass formula holds
\[
    \sum_{i = 1}^6 \frac{1}{|{\OO(W_i)}|} = \frac{505121}{12340763622899712000} = \mass(\cW_X). \qedhere
\]
\endproof

\begin{theorem} \label{thm:Barth-Peters}
If \(X\) is a K3 surface with transcendental lattice \(T_X \cong \bU \oplus \bU(2)\), then
\[
    |\cJ_X/{\Aut(X)}| = 7.
\]
\end{theorem}
\proof
By~\autoref{lem:O.hdg.Barth_Peters}, the subgroup \(\OO_\hdg^\sharp(T_X)\) is trivial, so by~\autoref{thm:main_multiplicity} the multiplicity of a frame~\(W \in \cW_X\) is equal to the index of \(\OO^\sharp(W)\) in \(\OO(T_X^\sharp)\). 
A finite set of generators of \(\OO(W)\) can be computed with the command \verb+orthogonal_group+ of the Sage class \verb+QuadraticForm+. Their images generate the subgroups \(\OO^\sharp(W)\). 
\endproof 

{\small
\begin{table}[t]
    \centering
    \caption{Lattices in the frame genus \(\cW_X\) of a K3 surface \(X\) with transcendental lattice \(T_X \cong \bU \oplus \bU(2)\).}
    \label{tab:genus_Barth_Peters}
    \begin{tabular}{lllllll}
    \toprule 
    \(W\) & \(W_\rootlattice\) & \(W/W_\rootlattice\) & \(|\Delta(W)|\) & \(|{\OO(W)}|\) & \(|{\fr_X^{-1}(W)}|\) \\
    \midrule
    \(W_1\) & \(\bD_8\bE_8\) & \(0\)                    & \(352\) & \(7191587192832000\)    & \(1\) \\
    \(W_2\) & \(\bD_{16}\)  & \(0\)                     & \(480\) & \(1371195958099968000\) & \(1\) \\
    \midrule 
    \(W_3\) & \(\bD_8^2\)   & \(\IZ/2\IZ\)              & \(224\) & \(53271016243200\)      & \(2\) \\
    \(W_4\) & \(\bA_1^2\bE_7^2\) & \(\IZ/2\IZ\)         & \(256\) & \(134842259865600\)     & \(1\) \\
    \(W_5\) & \(\bD_4 \bD_{12}\) & \(\IZ/2\IZ\)         & \(288\) & \(376702186291200\)     & \(1\) \\
    \midrule
    \(W_6\) & \(\bA_{15}\)  & \(\IZ \oplus (\IZ/2\IZ)\) & \(240\) & \(83691159552000\)      & \(1\) \\
    \bottomrule
    \end{tabular}
\end{table}
}

\subsection{Kummer surfaces associated to a product of elliptic curves} \label{sec:Oguiso}

Let now \(X\) denote a K3 surface with transcendental lattice 
\[
    T_X \cong \bU(2)^2.
\]

\begin{remark} 
By results of Nikulin \cite[Rmk.~2]{Nikulin:Kummer.surfaces} and Morrison \cite[Prop.~4.3(i)]{Morrison:K3.large.picard}, \(X\) is the Kummer surface associated to an abelian surface \(A\) with transcendental lattice \(T_A \cong \bU^2\). Consequently, the Néron--Severi lattice of \(A\) is isomorphic to \(\bU\), so \(A\) is isomorphic to the product of two non-isogenous elliptic curves, according to Ruppert's criterion~\cite{Ruppert:when}. 
Therefore, the number of orbits \(|\cJ_X/{\Aut(X)}|\) has  already been computed geometrically by Oguiso~\cite{Oguiso:jacobian.fibrations.Kummer}. 
Here we verify his results using our algebraic approach.
\end{remark}

We fix a basis \(t_1,\ldots,t_4 \in T_X\) so that the corresponding Gram matrix is
\[
    T \coloneqq {\small \begin{pmatrix} 0 & 2 & 0 & 0 \\ 2 & 0 & 0 & 0 \\ 0 & 0 & 0 & 2 \\ 0 & 0 & 2 & 0 \end{pmatrix}}.
\]
We make the following identification:
\begin{equation} \label{eq:O(TX).Oguiso}
    \OO(T_X) \cong \Set{A \in \GL_4(\IZ)}{ ATA^\intercal = T}.
\end{equation}

The classes modulo \(T_X\) of \(\frac12 t_1,\ldots ,\frac 12 t_4\) form a basis over \(\IF_2\) for \(T_X^\sharp\).
A computation shows that the group~\(\OO(T_X^\sharp)\) contains \(72\) elements and can be identified with the following subgroup of \(\GL_4(\IF_2)\):
\begin{equation} \label{eq:group.G.Oguiso} 
\OO(T_X^\sharp) \cong \left\langle {\small \begin{pmatrix} 1 & 0 & 0 & 0 \\
 0 & 1 & 0 & 0 \\
 0 & 0 & 0 & 1 \\
 0 & 0 & 1 & 0
\end{pmatrix}}, \,
{\small \begin{pmatrix} 0 & 1 & 0 & 1 \\
 1 & 0 & 0 & 0 \\
 0 & 0 & 0 & 1 \\
 1 & 0 & 1 & 0 
\end{pmatrix}} \right\rangle.
\end{equation}

Under identifications~\eqref{eq:O(TX).Oguiso} and \eqref{eq:group.G.Oguiso}, the natural homomorphism \(\OO(T_X) \rightarrow \OO(T_X^\sharp)\) is given by 
\(A \mapsto (A \mod 2)\).
We now consider its restriction to the subgroup of Hodge isometries \(\OO_\hdg(T_X)\).

\begin{lemma} \label{lem:ker=pm.id.Oguiso}
The kernel of \(\OO_\hdg(T_X) \rightarrow {\OO(T_X^\sharp)}\) is equal to \(\set{\pm \id}\).
\end{lemma}
\proof
Clearly, \(\set{\pm \id}\) is contained in the kernel. By \autoref{lem:ker.O_hdg} and \cite[Thm.~6.1]{Kondo:automorphisms.K3.which.act.trivially.on.Pic}, it suffices to show that there is no automorphism of \(X\) of order \(4\) acting trivially on \(S_X\). This follows from \cite[Thm.~1]{Schuett:K3.non-sp.aut.2-power}.
\endproof 

A computation with the GAP command \verb+ConjugacyClasses+ shows that the group \(\OO(T_X^\sharp)\) contains exactly \(9\) conjugacy classes. 
Let us define the following elements of \(\OO(T_X^\sharp)\):
\[
    h_1 \coloneqq {\small \begin{pmatrix}
     1 & 0 & 0 & 0 \\
     0 & 1 & 0 & 0 \\
     0 & 0 & 1 & 0 \\
     0 & 0 & 0 & 1
    \end{pmatrix}}, \quad
    h_2 \coloneqq {\small \begin{pmatrix}
     0 & 0 & 1 & 0 \\
     0 & 0 & 0 & 1 \\
     1 & 0 & 0 & 0 \\
     0 & 1 & 0 & 0
    \end{pmatrix}}, \quad
    h_3 \coloneqq {\small \begin{pmatrix}
     0 & 0 & 1 & 0 \\
     0 & 1 & 0 & 1 \\
     1 & 0 & 1 & 0 \\
     0 & 1 & 0 & 0 \\
    \end{pmatrix}}, \quad
    h_6 \coloneqq {\small \begin{pmatrix}
     0 & 0 & 0 & 1 \\
     1 & 1 & 1 & 1 \\
     0 & 1 & 0 & 1 \\
     1 & 0 & 0 & 1
    \end{pmatrix}}.
\]

\begin{lemma} \label{lem:order.2.Oguiso}
An element \(h \in {\OO_\hdg^\sharp(T_X)}\) of order \(2\) belongs to the conjugacy class of \(h_2\).
\end{lemma}
\proof
By~\autoref{prop:HLOY.Appendix.B} and \autoref{lem:ker=pm.id.Oguiso} we can suppose that \(h\) is the image of an element \(\eta \in \OO_\hdg(T_X)\) represented by a matrix \(A \in \GL_4(\IZ)\) of order \(4\) such that
\[ ATA^\intercal = T, \qquad A^2 = -I. \]

Since \(A\) has finite order, it is diagonalizable over \(\IC\)  with eigenvalues \(\ri = \sqrt{-1}\) and \(-\ri\).

Let \(V \subset T_X \otimes \IC\) be the eigenspace of \(A\) associated with \(\ri\).
Then \(V\) has dimension \(2\) and its complex conjugate \(\bar V\) is the eigenspace of \(A\) associated with \(-\ri\). 
As \(\OO_\hdg(T_X)\) acts faithfully on \(\IC \sigma_X\), we can assume that \(\sigma \coloneqq \sigma_X \in V\). 
Let us write 
\[
    \sigma = \sigma_1 t_1 + \sigma_2 t_2 + \sigma_3t_3 + \sigma_4t_4, \quad \sigma_i \in \IC,
\]
Certainly \(\sigma_1 \neq 0\) otherwise \(t_2 \in \sigma_X^\perp = S_X\), so we are free to rescale \(\sigma_1 = 1\). Moreover, from the relations \(\sigma^2 = 0\) and \(\sigma \cdot \bar \sigma > 0\) we obtain
\[
2(\sigma_1\sigma_2 + \sigma_3\sigma_4) = 0, \qquad \sigma_1\bar \sigma_2 + \sigma_2 \bar \sigma_1 + \sigma_3 \bar \sigma_4 + \sigma_4 \bar \sigma_3 > 0;
\]
substituting \(\sigma_1 = 1 \) and \(\sigma_2 = -\sigma_3\sigma_4\) we get \(\Im(\sigma_3)\Im(\sigma_4)>0\) (in particular, \(\Im(\sigma_3) \neq 0\)).

Let 
\[
    \tau = \tau_1 t_1 + \tau_2 t_2 + \tau_3 t_3 + \tau_4 t_4 \in V, \quad \tau_i \in \IC,
\]
be an eigenvector of \(A\) linearly independent from \(\sigma\).
Without loss of generality, we can substitute \(\tau\) with a linear combination of \(\sigma\) and \(\tau\) and suppose that \(\tau_1 = 0\) 

Note that \(V\) is a totally isotropic subspace of \(T_X \otimes \IC\), as
\[
    x \cdot y = \eta(x) \cdot \eta(y) = (\ri x) \cdot (\ri y) = - x \cdot y, \quad \text{for all \(x,y \in V\)}.
\]
From \(\tau^2 = 0\) it follows that \(2\tau_3\tau_4 = 0\). Up to exchanging \(t_3\) with \(t_4\) we can suppose that \(\tau_3 = 0\). 
From the relation \(\sigma \cdot \tau = 0\) we obtain \(\tau_2 = -\sigma_3 \tau_4\). Again without loss of generality, we can rescale \(\tau_4 = 1\).

Imposing that \(\sigma,\tau \in V\) and \(\bar\sigma,\bar\tau \in \bar V\) and recalling that \(A\) has integer coefficients, we infer with elementary but rather tedious computations that there exist \(m_1,m_2 \in \IZ\) such that
\[
    A = {\small \begin{pmatrix}
     m_1 & 0 & -(m_1^2+1)/m_2 & 0 \\
     0 & -m_1 & 0 & -m_2 \\
     m_2 & 0 & -m_1 & 0 \\
     0 & (m_1^2+1)/m_2 & 0 & m_1
    \end{pmatrix}}, \quad m_1 \coloneqq -\frac{\Re(\sigma_3)}{\Im(\sigma_3)}, \quad m_2 \coloneqq \frac{1}{\Im(\sigma_3)}.
\]
Note that \(m_1,m_2\) cannot be both even because \(h = (A \mod 2)\) is an invertible matrix. Looking at the possible parities of \(m_1,m_2\), we see that \(h\) is equal to one of the following elements:
\[
    h_2 \coloneqq {\small \begin{pmatrix}
    0 & 0 & 1 & 0 \\
    0 & 0 & 0 & 1 \\
    1 & 0 & 0 & 0 \\
    0 & 1 & 0 & 0
    \end{pmatrix}}, \quad
    {\small \begin{pmatrix}
     1 & 0 & 0 & 0 \\
     0 & 1 & 0 & 1 \\
     1 & 0 & 1 & 0 \\
     0 & 0 & 0 & 1
    \end{pmatrix}}, \quad
    {\small \begin{pmatrix}
     1 & 0 & 1 & 0 \\
     0 & 1 & 0 & 0 \\
     0 & 0 & 1 & 0 \\
     0 & 1 & 0 & 1
    \end{pmatrix}}, 
\]
which all belong to the same conjugacy class.
\endproof 

\begin{lemma} \label{lem:order.3.Oguiso}
An element \(h \in \OO_\hdg^\sharp(T_X)\) of order \(3\) belongs to the conjugacy class of \(h_3\).
\end{lemma}
\proof
By~\autoref{lem:ker=pm.id.Oguiso}, we can suppose that \(h\) is the image of an element in \(\OO_\hdg(T_X)\) represented by a matrix \(A \in \GL_4(\IZ)\) such that
\[ ATA^\intercal = T, \qquad A^3 = I. \]
Since \(A\) has finite order \(3\), it is diagonalizable over \(\IC\)  with eigenvalues \(1,\omega,\bar \omega\), where \(\omega\) denotes a third root of unity. 
Let \(\zeta\) be a primitve \(12\)th root of unity such that \(\zeta^4 = \omega\) and \(\zeta^3 = \ri\).

Let \(W \subset T_X \otimes \IC\) be the eigenspace of \(A\) associated with \(\omega\). We can suppose that \(\sigma \coloneqq \sigma_X \in W\).  
We must differentiate two cases: either \(\dim(W) = 1\) or \(\dim(W) = 2\).

Suppose first that \(\dim(W) = 2\), let \(\tau \in W\) be linearly independent from \(\sigma\).
Without loss of generality we can suppose as in the proof \autoref{lem:order.2.Oguiso} that \(\sigma_1 = 1\), \(\sigma_2 = -\sigma_3\sigma_4\), \(\tau_1 = \tau_3 = 0\), \(\tau_2 = -\sigma_3\tau_4\), \(\tau_4 = 1\).

Imposing that \(\sigma,\tau \in W\) and \(\bar\sigma,\bar\tau \in \bar W\) we find (again after some tedious computations) that
\[
A = {\small \begin{pmatrix}
m_1 & 0 & -(m_1^2 + m_1 +1)/m_2 & 0 \\
0 & -m_1-1 & 0 & -m_2 \\ 
m_2 & 0 & -m_1 -1 & 0 \\
0 & (m_1^2 + m_1 +1)/m_2 & 0 & m_1 \\
\end{pmatrix}},
\]
where
\[
m_1 \coloneqq \frac{(\zeta^3-2\zeta)\Re(\sigma_3) - 4\Im(\sigma_3)}{4\Im(\sigma_3)}, \quad m_2 \coloneqq \frac{2\zeta-\zeta^3}{2\Im(\sigma_3)}.
\]
Looking at the possible parities of \(m_1,m_2\), we see that \(h\) is equal to one of the following elements:
\[
    h_3 \coloneqq {\small \begin{pmatrix}
    0 & 0 & 1 & 0 \\
    0 & 1 & 0 & 1 \\
    1 & 0 & 1 & 0 \\
    0 & 1 & 0 & 0
    \end{pmatrix}}, \quad 
    {\small \begin{pmatrix}
     1 & 0 & 1 & 0 \\
     0 & 0 & 0 & 1 \\
     1 & 0 & 0 & 0 \\
     0 & 1 & 0 & 1
    \end{pmatrix}}, 
\]
which are conjugate to each other.

In the case that \(\dim(W) = 1\) we can assume that \(W\) is generated by \(\sigma\) and that the eigenspace \(U\) associated with \(1\) is generated by \(\tau\) and \(\upsilon\). 
We note that \(W\) and \(U\) are orthogonal to each other, as 
\[
    w \cdot u = \eta(w) \cdot \eta(u) = \omega w \cdot u, \quad \text{for all \(u \in U, w \in W\)}.
\]
(Note, though, that \(U\) is not necessarily a totally isotropic subspace.) With similar computations we see that the matrix \(A\) is forced to have non-integral coefficients. Thus, this case is impossible.
\endproof 

\begin{lemma} \label{lem:order.4.Oguiso}
The group \(\OO_\hdg^\sharp(T_X)\) contains no element of order \(4\).
\end{lemma}
\proof 
By inspection of the conjugacy classes of \(\OO(T_X^\sharp)\), we see that there is only one class containing elements of order~\(4\), namely the class of
\[
    h_4 \coloneqq {\small \begin{pmatrix}
     0 & 0 & 0 & 1 \\
     0 & 0 & 1 & 0 \\
     1 & 0 & 0 & 0 \\
     0 & 1 & 0 & 0
    \end{pmatrix}}.
\]
We see, though, that \(h_4^2\) is not conjugate to \(h_2\), so we conclude by \autoref{lem:order.2.Oguiso}.
\endproof 

\begin{proposition} \label{prop:Hodge.isometries.Oguiso}
The subgroup \({\OO_\hdg^\sharp(T_X)}\) is a cyclic group of order \(1,2,3\) or \(6\), in which case it is generated by a conjugate of \(h_1,h_2,h_3\) or \(h_6\), respectively.
\end{proposition}
\proof
Let \(h\) be a generator of \({\OO_\hdg^\sharp(T_X)}\) and let \(m\) be its order.

By \autoref{prop:HLOY.Appendix.B} it holds \(\varphi(|{\OO_\hdg(T_X)}|) \mid \rank(T_X) = 4\), therefore
\(
    |{\OO_\hdg(T_X)}| \in \set{2,4,6,8,10,12}.
\)
It follows from \autoref{lem:ker=pm.id.Oguiso} that \(m = |{\OO_\hdg(T_X)}|/2 \in \set{1,2,3,4,5,6}\).

We can exclude \(m = 5\) because \(5 \nmid |{\OO(T_X^\sharp)}| = 72\) and we can exclude \(m = 4\) because of \autoref{lem:order.4.Oguiso}. 
All in all, \(m \in \set{1,2,3,6}\) as claimed.

If \(m = 1\), then obviously \(h = h_1\).
If \(m = 2\), then \(h\) is conjugate to \(h_2\) by \autoref{lem:order.2.Oguiso}. 
If \(m = 3\), then \(h\) is conjugate to \(h_3\) by \autoref{lem:order.3.Oguiso}. 
Finally, if \(m = 6\), then by inspection of the conjugacy classes we see that \(h\) is conjugate to either \(h_6\) or
\[
    h_6' \coloneqq {\small \begin{pmatrix}
     0 & 0 & 0 & 1 \\
     0 & 1 & 1 & 0 \\
     1 & 0 & 0 & 1 \\
     0 & 1 & 0 & 0 \\
    \end{pmatrix}}.
\]
Thanks to \autoref{lem:order.2.Oguiso}, we can exclude \(h_6'\), because \((h_6')^3\) is not conjugate to \(h_2\).
\endproof

\begin{proposition}
The frame genus \(\cW_X\) contains exactly \(11\) isomorphism classes, listed in \autoref{tab:genus_Oguiso}, whose Gram matrices are contained in the arXiv ancillary file \verb+genus_Oguiso.sage+. 
\end{proposition}
\proof 
We apply the Kneser--Nishiyama method with \(T_0 = \bD_4^2\).
Thanks to \cite[Lemma 41(i) and Corollary 4.6(i)]{Nishiyama:Jacobian.fibrations.K3.MW} we find exatly \(10\) Niemeier lattices~\(N\) such that there exists a primitive embedding \(T_0 \hookrightarrow N_\rootlattice\).
All primitive embeddings are unique, except when \(N_\rootlattice = \bD_{10}\bE_7^{2}\), in which case there are two primitive embeddings (corresponding to \(W_5\) and \(W_8\) in \autoref{tab:genus_Oguiso}). 
Therefore, we obtain \(11\) different frames: \(\cW_X = \set{W_1,\ldots,W_{11}}\).
The list is complete because the mass formula  holds (see \autoref{rem:mass}):
\[
\sum_{i = 1}^{11} \frac{1}{|{\OO(W_i)}|} = \frac{64150367}{1708721117016883200} = \mass(\cW_X). \qedhere 
\]
\endproof

We introduce the following subgroups of \(G\):
\begin{align*} 
K_8 & \coloneqq \left\langle 
{\small \begin{pmatrix}
    1 & 0 & 0 & 1 \\ 0 & 1 & 0 & 0 \\ 0 & 1 & 1 & 0 \\ 0 & 0 & 0 & 1
\end{pmatrix}},
{\small \begin{pmatrix}
    1 & 0 & 0 & 0 \\ 0 & 1 & 0 & 0 \\ 0 & 0 & 0 & 1 \\ 0 & 0 & 1 & 0
\end{pmatrix}} 
\right\rangle, \\
K_{12} & \coloneqq \left\langle 
{\small \begin{pmatrix}
    1 & 1 & 1 & 1 \\ 1 & 0 & 0 & 0 \\ 1 & 0 & 0 & 1 \\ 1 & 0 & 1 & 0
\end{pmatrix}},
{\small \begin{pmatrix} 
    1 & 0 & 0 & 0 \\ 1 & 1 & 1 & 1 \\ 1 & 0 & 0 & 1 \\ 1 & 0 & 1 & 0
\end{pmatrix}},
{\small \begin{pmatrix}
    1 & 0 & 0 & 0 \\ 0 & 1 & 0 & 0 \\ 0 & 0 & 0 & 1 \\ 0 & 0 & 1 & 0        
\end{pmatrix}}
\right\rangle, \\            
K_{36} & \coloneqq \left\langle 
{\small \begin{pmatrix}
1 & 0 & 0 & 0 \\ 1 & 1 & 1 & 1 \\ 1 & 0 & 1 & 0 \\ 1 & 0 & 0 & 1
\end{pmatrix}},
{\small \begin{pmatrix}
    0 & 1 & 1 & 0 \\ 0 & 0 & 0 & 1 \\ 1 & 0 & 0 & 1 \\ 0 & 1 & 0 & 0
\end{pmatrix}},
{\small \begin{pmatrix} 
    1 & 0 & 0 & 1 \\ 0 & 1 & 0 & 0 \\ 0 & 1 & 1 & 0 \\ 0 & 0 & 0 & 1
\end{pmatrix}} 
    \right\rangle.
\end{align*}
The subgroups \(K_8,K_{12},K_{36}\) contain respectively \(8,12,36\) elements. Note that neither \(K_8\) nor \(K_{12}\) are normal subgroups of \(G\), but \(K_{36}\) is, as it has index \(2\).

\begin{proposition} \label{prop:O^sharp(W).Oguiso}
The subgroup \(\OO^\sharp(W) \subset \OO(T_X^\sharp)\) is conjugate 
to~\(K_8\) if \(W \in \set{ W_1,W_6,W_7,W_{11} }\), 
to~\(K_{12}\) if \(W \in \set{W_2,W_8,W_{10}}\), 
to~\(K_{36}\) if \(W = W_4\), 
and to~\(G\) if \(W \in \set{W_3,W_5,W_9}\).
\end{proposition}%
\proof 
A finite set of generators of \(\OO(W)\) can be computed with the command \verb+orthogonal_group+ of the Sage class \verb+QuadraticForm+. Their images generate the subgroups \(\OO^\sharp(W)\). 
\endproof 

We now have all ingredients to compute the multiplicities of the frames \(W \in \cW_X\).

\begin{theorem} \label{thm:Oguiso}
If \(X\) is a K3 surface with transcendental lattice \(T_X \cong \bU(2)^2\), then one of the following cases holds:
\[
    |\cJ_X/{\Aut(X)}| = \begin{cases} 59 & \text{if \(|{\OO^\sharp_\hdg(T_X)}| = 1\)}, \\ 38 & \text{if \(|{\OO^\sharp_\hdg(T_X)}| = 2\)}, \\ 23 & \text{if \(|{\OO^\sharp_\hdg(T_X)}| = 3\)}, \\ 16 & \text{if \(|{\OO^\sharp_\hdg(T_X)}| = 6\)}. \end{cases} 
\]
\end{theorem}
\proof
We apply \autoref{thm:main_multiplicity} to compute the multiplicities of the frames \(W \in \cW_X\) using either formula~\eqref{eq:Cauchy-Frobenius} or the GAP function \verb+DoubleCosets+.
Thanks to \autoref{prop:Hodge.isometries.Oguiso} and \autoref{prop:O^sharp(W).Oguiso}, the multiplicities are given by \(|H\backslash G /K|\), where \(G\) is the group defined in \eqref{eq:group.G.Oguiso}, \(H\) is the cyclic subgroup generated by \(h_1,h_2,h_3\) or \(h_6\) and \(K \in \set{K_8, K_{12}, K_{36}}\).
Our results are listed in \autoref{tab:multiplicities.Oguiso}.
\endproof%

{\small
\begin{table}
    \centering
    \caption{Lattices in the frame genus \(\cW_X\) of a K3 surface \(X\) with transcendental lattice \(T_X \cong \bU(2)^2\), numbered according to Oguiso (cf.~\cite[Table~2]{Oguiso:jacobian.fibrations.Kummer}).}
    \label{tab:genus_Oguiso}
    \begin{tabular}{lllllll}
    \toprule
    \(W\) & \(N_\rootlattice\) & \(W_\rootlattice\) & \(W/W_\rootlattice\) & \(|\Delta(W)|\) & \(|{\OO(W)}|\) & \(|{\OO^\sharp(W)}|\) \\
    \midrule 
    \(W_{11}\) & \(\bD_{12}^2\) & \(\bD_{8}^2\) & \(0\) & \(224\) & \(213084064972800\) & \(8\) \\
    \(W_{10}\) & \(\bD_{16}\bE_8\) & \(\bD_{4} \bD_{12}\) & \(0\) & \(288\) & \(2260213117747200\) & \(12\) \\
    \(W_9\) & \(\bE_8^3\) & \(\bD_4^2 \bE_8\) & \(0\) & \(288\) & \(1849265278156800\) & \(72\) \\
    \midrule 
    \(W_7\) & \(\bD_8^3\) & \(\bD_4^2 \bD_8\) & \(\IZ/2\IZ\) & \(160\) & \(1522029035520\) & \(8\) \\
    \(W_8\) & \(\bD_{10}\bE_7^2\) & \(\bA_1^3 \bD_6 \bE_7\) & \(\IZ/2\IZ\) & \(192\) & \(6421059993600\) & \(12\) \\
    \midrule
    \(W_4\) & \(\bD_4^6\) & \(\bD_4^4\) & \((\IZ/2\IZ)^2\) & \(96\) & \(195689447424\) & \(36\) \\
    \(W_6\) & \(\bD_6^4\) & \(\bA_1^4\bD_6^2\) & \((\IZ/2\IZ)^2\) & \(128\) & \(135895449600\) & \(8\) \\
    \(W_5\) & \(\bD_{10}\bE_7^2\) & \(\bA_1^6 \bD_{10}\) & \((\IZ/2\IZ)^2\) & \(192\) & \(8561413324800\) & \(72\) \\
    \midrule
    \(W_1\) & \(\bA_7^2\bD_5^2\) & \(\bA_7^2\) & \(\IZ^2 \oplus (\IZ/2\IZ)\) & \(112\) & \(52022476800\) & \(8\) \\
    \(W_2\) & \(\bA_{11}\bD_7\bE_6\) & \(\bA_3 \bA_{11} \) & \(\IZ^2 \oplus (\IZ/2\IZ)\) & \(144\) & \(275904921600\) & \(12\) \\
    \midrule 
    \(W_3\) & \(\bE_6^4\) & \(\bE_6^2\) & \(\IZ^4\) & \(144\) & \(773967052800\) & \(72\) \\
    \bottomrule
    \end{tabular}
\end{table}

\begin{table}
    \centering
    \caption{Multiplicities of the frames \(W \in \cW_X\) listed in \autoref{tab:genus_Oguiso}.}
    \label{tab:multiplicities.Oguiso}
    \begin{tabular}{c|ccccccccccc|c}
    \toprule
        \(|{\OO_\hdg^\sharp(T_X)}|\) & \(W_1\) & \(W_2\) & \(W_3\) & \(W_4\) & \(W_5\) & \(W_6\) & \(W_7\) & \(W_8\) & \(W_9\) & \(W_{10}\) & \(W_{11}\) & sum \\
    \midrule
        \(1\) & \(9\) & \(6\) & \(1\) & \(2\) & \(1\) & \(9\) & \(9\) & \(6\) & \(1\) & \(6\) & \(9\) & \(59\) \\
        \(2\) & \(6\) & \(3\) & \(1\) & \(2\) & \(1\) & \(6\) & \(6\) & \(3\) & \(1\) & \(3\) & \(6\)& \(38\) \\
        \(3\) & \(3\) & \(2\) & \(1\) & \(2\) & \(1\) & \(3\) & \(3\) & \(2\) & \(1\) & \(2\) & \(3\) & \(23\) \\
        \(6\) & \(2\) & \(1\) & \(1\) & \(2\) & \(1\) & \(2\) & \(2\) & \(1\) & \(1\) & \(1\) & \(2\) & \(16\) \\
    \bottomrule
    \end{tabular}
\end{table}
}
\begin{remark} \label{rem:geom.interpr.Oguiso}
With \autoref{thm:Oguiso} we confirm independently Oguiso's results \cite{Oguiso:jacobian.fibrations.Kummer}.
Following his notation, we denote by \(E_\tau\) the elliptic curve with period \(\tau\).
Comparing \autoref{tab:multiplicities.Oguiso} with \cite[Table B]{Oguiso:jacobian.fibrations.Kummer}, we find a geometrical interpretation of the order of \({\OO_\hdg^\sharp(T_X)}\). Indeed, it holds
\[
    |{\OO_\hdg^\sharp(T_X)}| = \begin{cases}
    1 & \text{if \(X \cong \Km(E_{\rho} \times E_{\rho'})\)}, \\
    2 & \text{if \(X \cong \Km(E_{\sqrt{-1}} \times E_{\rho'})\)}, \\
    3 & \text{if \(X \cong \Km(E_\rho \times E_{(-1+\sqrt{-3})/2})\)}, \\
    6 & \text{if \(X \cong \Km(E_{\sqrt{-1}} \times E_{(-1+\sqrt{-3})/2})\)},
    \end{cases}
\]
where \(E_\rho,E_{\rho'}\) are non-isogenous elliptic curves without complex multiplication.
\end{remark}

\subsection{Jacobian Kummer surfaces} \label{sec:Kumar}

Here let \(X\) denote a K3 surface with transcendental lattice 
\[
    T_X \cong \bU(2)^2 \oplus [-4].
\]

\begin{remark}
By results of Kumar~\cite{Kumar:elliptic.fibrations},
the Kummer surface \(X\) of the Jacobian of a curve \(C\) of genus \(2\) without extra endomorphisms satisfies \(T_X \cong \bU(2)^2 \oplus [-4]\).
In the same paper, he classified the elliptic fibrations up to the action of \(\Aut(X)\) and the permutations of the Weierstrass points of~\(C\), establishing that \(|\cW_X| = 25\) and finding all pairs \((W_\rootlattice,W/W_\rootlattice), W \in \cW_X\) (cf. \cite[Thm.~2]{Kumar:elliptic.fibrations}).
Before Kumar's work, some special elliptic fibrations had been found by Keum~\cite{Keum:two.extremal,Keum:erratum:two.extremal}, Shioda~\cite{Shioda:classical.Kummer} and Kumar himself~\cite{Kumar:K3.associated}.
\end{remark}

\begin{proposition} \label{prop:genus.Kumar}
The frame genus \(\cW_X\) contains exactly \(25\) isomorphism classes, listed in \autoref{tab:genus_Kumar}, whose Gram matrices are contained in the arXiv ancillary file \verb+genus_Kumar.sage+. 
\end{proposition}
\proof
We computed \(\cW_X\) using Kneser's neighbor method (\autoref{rem:Kneser.neighbor}).
\endproof 

\begin{remark}
Alternatively, one could compute \(\cW_X\) applying the Kneser--Nishiyama method with \(T_0 = \bD_4 \oplus \bD_4 \oplus [-4]\). This lattice, though, is not generated by its roots, so one cannot use Nishiyama's results \cite{Nishiyama:Jacobian.fibrations.K3.MW}.
\end{remark}

\begin{theorem} \label{thm:Kumar}
If \(X\) is a K3 surface with transcendental lattice \(T_X \cong \bU(2)^2 \oplus [-4]\), then 
\[
    |\cJ_X/{\Aut(X)}| = 491.
\]
\end{theorem}
\proof
It holds \(|{\OO(T_X^\sharp)}| = 1440\). Since \(\rank(T_X)\) is odd, \(|{\fr_X^{-1}(W)}|\) is equal to the index of \(\OO^\sharp(W)\) in \(\OO(T_X^\sharp)\) (\autoref{rmk:rk.TX.odd}).
A finite set of generators of \(\OO(W)\) can be computed using the command \verb+orthogonal_group+ of the Sage class \verb+QuadraticForm+. Their images generate \(\OO^\sharp(W)\).
\endproof 
{\small 
\begin{table}
    \centering
    \caption{Lattices in the frame genus \(\cW_X\) of a K3 surface \(X\) with transcendental lattice \(T_X \cong \bU(2)^2 \oplus [-4]\), numbered according to Kumar (cf.~\cite[§3.3]{Kumar:elliptic.fibrations}).}
    \label{tab:genus_Kumar}
    \begin{tabular}{llllll}
    \toprule
    \(W\) & \(W_\rootlattice\) & \(W/W_\rootlattice\) & \(|\Delta(W)|\) & \(|{\OO(W)}|\) & \(|{\fr_X^{-1}(W)}|\) \\
    \midrule 
    \(W_{17}\) & \(\bD_4^2\bD_7\)  & \(0\) & \(132\) & \(1712282664960\) & \(10\) \\
    \(W_{15}\) & \(\bA_3\bD_4\bD_8\) & \(0\) & \(148\) & \(570760888320\) & \(60\) \\
    \midrule
    \(W_4\) & \(\bA_3\bD_4^3\) & \(\IZ/2\IZ\) & \(84\) & \(16307453952\) & \(15\) \\
    \(W_9\) & \(\bA_1^4\bD_5\bD_6\) & \(\IZ/2\IZ\) & \(108\) & \(67947724800\) & \(15\) \\
    \(W_{18}\) & \(\bA_1^5\bA_3\bE_7\) & \(\IZ/2\IZ\) & \(148\) & \(535088332800\) & \(6\) \\
    \(W_{23}\) & \(\bA_1^6\bD_9\) & \(\IZ/2\IZ\) & \(156\) & \(8561413324800\) & \(1\) \\
    \midrule
    \(W_3\) & \(\bA_1^6\bA_3\bD_6\) & \((\IZ/2\IZ)^2\) & \(84\) & \(5096079360\) & \(10\) \\
    \midrule
    \(W_8\) & \(\bA_1^2\bD_6^2\) & \(\IZ\) & \(124\) & \(67947724800\) & \(45\) \\
    \(W_{16}\) & \(\bA_1^3 \bD_4 \bE_7\) & \(\IZ\) & \(156\) & \(321052999680\) & \(20\) \\
    \(W_{20}\) & \(\bA_1^4 \bD_{10}\) & \(\IZ\) & \(188\) & \(2853804441600\) & \(15\) \\
    \(W_{24}\) & \(\bA_1^6\bE_8\) & \(\IZ\) & \(252\) & \(64210599936000\) & \(1\) \\
    \midrule
    \(W_2\) & \(\bA_1^4\bD_4\bD_6\) & \(\IZ \oplus (\IZ/2\IZ)\) & \(92\) & \(1698693120\) & \(60\) \\
    \(W_{11}\) & \(\bA_1^6\bD_8\) & \(\IZ \oplus (\IZ/2\IZ)\) & \(124\) & \(31708938240\) & \(15\) \\
    \midrule
    \(W_1\) & \(\bA_1^6\bD_4^2\) & \(\IZ \oplus (\IZ/2\IZ)^2\) & \(60\) & \(452984832\) & \(15\) \\
    \midrule
    \(W_5\) & \(\bA_3^2\bA_7\) & \(\IZ^2 \oplus (\IZ/2\IZ)\) & \(80\) & \(743178240\) & \(45\) \\
    \midrule
    \(W_7\) & \(\bA_3^4\) & \(\IZ^3 \oplus (\IZ/2\IZ)\) & \(48\) & \(127401984\) & \(15\) \\
    \(W_{10}\) & \(\bA_1^2 \bA_5^2\) & \(\IZ^3 \oplus (\IZ/2\IZ)\) & \(64\) & \(99532800\) & \(60\) \\
    \(W_6\) & \(\bA_7\bA_3\bA_1^2\) & \(\IZ^3 \oplus (\IZ/2\IZ)\) & \(72\) & \(743178240\) & \(15\) \\
    \(W_{21}\) & \(\bA_1^3\bA_9\) & \(\IZ^3 \oplus (\IZ/2\IZ)\) & \(96\) & \(2090188800\) & \(20\) \\
    \midrule
    \(W_{12}\) & \(\bD_5 \bE_6\) & \(\IZ^4\) & \(112\) & \(9555148800\) & \(15\) \\
    \midrule
    \(W_{14}\) & \(\bA_5^2\) & \(\IZ^5\) & \(60\) & \(298598400\) & \(10\) \\
    \(W_{19}\) & \(\bD_5^2\) & \(\IZ^5\) & \(80\) & \(707788800\) & \(15\) \\
    \(W_{22}\) & \(\bA_1\bA_9\) & \(\IZ^5\) & \(92\) & \(3483648000\) & \(6\) \\
    \(W_{13}\) & \(\bD_4\bE_6\) & \(\IZ^5\) & \(96\) & \(28665446400\) & \(1\) \\
    \midrule
    \(W_{25}\) & \(\bD_9\) & \(\IZ^6\) & \(144\) & \(133772083200\) & \(1\) \\
    \bottomrule
    \end{tabular}
\end{table}
}

\subsection{Double covers of the plane ramified over six lines} \label{sec:Kloosterman}

Here let \(X\) denote a K3 surface with transcendental lattice 
\[
    T_X \cong \bU(2)^2 \oplus [-2]^2.
\]

\begin{remark}
By results of Matsumoto, Sasaki and Yoshida~\cite{Matsumoto.Sasaki.Yoshida:monodromy.4-parameter}, \(X\) is the minimal resolution of the double cover of \(\IP^2\) ramified over \(6\) lines, and,
conversely, such a resolution satisfies \(T_X \cong \bU(2)^2 \oplus [-2]^2\), provided that \(\rank(S_X) = 16\).
Kloosterman \cite{Kloosterman:classification.jac.ell.fibr} classified the fiber types and the Mordell--Weil groups of all jacobian elliptic fibrations on \(X\).
Here we refine Kloosterman's classification and compute the number of orbits \(|\cJ_X/{\Aut(X)}|\).
\end{remark}

We fix a basis \(t_1,\ldots,t_6 \in T_X\) so that the corresponding Gram matrix is
\[
    T \coloneqq {\small \begin{pmatrix} 0 & 2 & 0 & 0 & 0 & 0 \\ 2 & 0 & 0 & 0 & 0 & 0 \\ 0 & 0 & 0 & 2 & 0 & 0 \\ 0 & 0 & 2 & 0 & 0 & 0 \\ 0 & 0 & 0 & 0 & -2 & 0 \\ 0 & 0 & 0 & 0 & 0 & -2 \end{pmatrix}}.
\]
We make the following identification:
\begin{equation} \label{eq:O(TX).Kloosterman}
    \OO(T_X) \cong \Set{A \in \GL_6(\IZ)}{ ATA^\intercal = T}.
\end{equation}

The classes modulo \(T_X\) of \(\frac12 t_1,\ldots \frac 12 t_6\) form a basis over \(\IF_2\) for \(T_X^\sharp\).
A computation shows that the group~\(\OO(T_X^\sharp)\) contains \(1440 = 2^5 \cdot 3^2 \cdot 5\) elements and can be identified with the following subgroup of \(\GL_6(\IF_2)\):
\begin{equation} \label{eq:group.G.Kloosterman}
\OO(T_X^\sharp) \cong \left\langle 
{\small \begin{pmatrix} 0 & 0 & 1 & 1 & 1 & 1 \\ 1 & 0 & 1 & 0 & 0 & 0 \\ 1 & 1 & 1 & 1 & 0 & 0 \\ 0 & 1 & 1 & 1 & 1 & 1 \\ 0 & 1 & 0 & 1 & 1 & 0 \\ 0 & 1 & 0 & 1 & 0 & 1 \end{pmatrix}}, \,
{\small \begin{pmatrix} 0 & 0 & 0 & 1 & 0 & 0 \\ 0 & 0 & 1 & 0 & 0 & 0 \\ 1 & 0 & 0 & 0 & 0 & 0 \\ 0 & 1 & 0 & 0 & 0 & 0 \\ 0 & 0 & 0 & 0 & 1 & 0 \\ 0 & 0 & 0 & 0 & 0 & 1 \end{pmatrix}}, \,
{\small \begin{pmatrix} 1 & 0 & 0 & 0 & 0 & 0 \\ 0 & 1 & 0 & 0 & 0 & 0 \\ 0 & 0 & 1 & 0 & 0 & 0 \\ 0 & 0 & 0 & 1 & 0 & 0 \\ 0 & 0 & 0 & 0 & 0 & 1 \\ 0 & 0 & 0 & 0 & 1 & 0 \end{pmatrix}} \right\rangle.
\end{equation}

Under identifications~\eqref{eq:O(TX).Kloosterman} and \eqref{eq:group.G.Kloosterman}, the natural homomorphism \(\OO(T_X) \rightarrow \OO(T_X^\sharp)\) is given by \(A \mapsto A \mod 2\).

\begin{lemma} \label{lem:ker=pm.id.Kloosterman}
The kernel of \(\OO_\hdg(T_X) \rightarrow \OO(T_X^\sharp)\) is equal to \(\set{ \pm \id }\).
\end{lemma}
\proof
Clearly, \(\set{\pm \id}\) is contained in the kernel. By \autoref{lem:ker.O_hdg} and \cite[Thm.~6.1]{Kondo:automorphisms.K3.which.act.trivially.on.Pic}, it suffices to show that there is no automorphism of \(X\) of order \(4\) acting trivially on \(S_X\). This follows from \cite[Main theorem on p.~18]{Taki:classification.non-sp.aut.K3.which.act.trivially}, noticing that in our case it holds \(\delta_{S_X}=1\) (cf. \cite[Definition~2.2]{Taki:classification.non-sp.aut.K3.which.act.trivially}), as \((t_6/2)^2 = -1/2\).
\endproof

A computation with the GAP command \verb+ConjugacyClasses+ shows that the group \(\OO(T_X^\sharp)\) contains exactly \(22\) conjugacy classes. 
Let us define the following elements of \(\OO(T_X^\sharp)\):
\[
    h_1 \coloneqq {\small \begin{pmatrix}
     1 & 0 & 0 & 0 & 0 & 0 \\
     0 & 1 & 0 & 0 & 0 & 0 \\
     0 & 0 & 1 & 0 & 0 & 0 \\
     0 & 0 & 0 & 1 & 0 & 0 \\
     0 & 0 & 0 & 0 & 1 & 0 \\
     0 & 0 & 0 & 0 & 0 & 1
    \end{pmatrix}}, \quad
    h_2 \coloneqq {\small \begin{pmatrix}
     0 & 0 & 0 & 1 & 0 & 0 \\
     0 & 0 & 1 & 1 & 1 & 1 \\
     1 & 1 & 0 & 0 & 1 & 1 \\
     1 & 0 & 0 & 0 & 0 & 0 \\
     1 & 0 & 0 & 1 & 0 & 1 \\
     1 & 0 & 0 & 1 & 1 & 0
    \end{pmatrix}}, \quad
    h_3 \coloneqq {\small \begin{pmatrix}
     0 & 0 & 1 & 0 & 0 & 0 \\
     0 & 1 & 0 & 1 & 0 & 0 \\
     1 & 0 & 1 & 0 & 0 & 0 \\
     0 & 1 & 0 & 0 & 0 & 0 \\
     0 & 0 & 0 & 0 & 1 & 0 \\
     0 & 0 & 0 & 0 & 0 & 1
    \end{pmatrix}}.
\]

\begin{lemma} \label{lem:order.2.Kloosterman}
An element \(h \in \OO_\hdg^\sharp(T_X)\) of order \(2\) belongs to the conjugacy class of \(h_2\).
\end{lemma}
\proof 
As in \autoref{lem:order.2.Oguiso}, thanks to \autoref{lem:ker=pm.id.Kloosterman} we can suppose that \(h = (A \mod 2)\), where \(A \in \GL_6(\IZ)\) satisfies \(ATA^\intercal = T\) and \(A^2 = -I\).

Let \(V \subset T_X \otimes \IC\) be the eigenspace of \(A\) associated with \(\ri\). Necessarily, \(\dim(V) = 3\) and we can suppose that \(\sigma \coloneqq \sigma_X \in V\). 
Let \(\sigma,\tau,\upsilon\) be a basis of \(V\) over \(\IC\). 
We write
\begin{align*}
    \sigma & = \sigma_1 t_1 + \sigma_2 t_2 + \sigma_3 t_3 + \sigma_4 t_4 + \sigma_5 t_5 + \sigma_6 t_6, \quad \sigma_i \in \IC, \\  
    \tau & = \tau_1 t_1 + \tau_2 t_2 + \tau_3 t_3 + \tau_4 t_4 + \tau_5 t_5 + \tau_6 t_6, \quad \tau_i \in \IC, \\  
    \upsilon & = \upsilon_1 t_1 + \upsilon_2 t_2 + \upsilon_3 t_3 + \upsilon_4 t_4 + \upsilon_5 t_5 + \upsilon_6 t_6, \quad \upsilon_i \in \IC. 
\end{align*}
Again, \(V\) is a totally isotropic subspace with respect to the bilinear form induced by \(T_X\). Therefore, upon substituting \(t_3\) with \(t_4\), or \(\tau\) and \(\upsilon\) with linear combinations of \(\sigma,\tau,\upsilon\), we can make the following assumptions:
\begin{align*}
\sigma_1 &= 1, & 
\upsilon_1 &= 0 & 
\sigma_2 & = \sigma_3\sigma_4 + \tfrac12 \sigma_5^2 + \tfrac 12\sigma_6^2, & 
\upsilon_2 &= -\sigma_4\upsilon_3 - \sigma_3\upsilon_4 + \sigma_5\upsilon_5 + \sigma_6\upsilon_6, \\
\tau_1 &= 0, & 
\tau_5 &= 0, & 
\tau_2 &= -\sigma_4\tau_3 - \sigma_3\tau_4 + \sigma_6\tau_6, &
\upsilon_6 &= \tau_6\upsilon_3 + \ri \upsilon_5.
\end{align*}

Imposing that \(\sigma,\tau,\upsilon \in V\) and \(\bar\sigma,\bar\tau,\bar\upsilon \in \bar V\) and recalling that \(A\) has integer coefficients, we get with
elementary but rather tedious computations that there exist \(m_1,\ldots,m_6 \in \IZ\) such that
\[
A = {\small \begin{pmatrix} 
m_1 & 0 & * & * & * & * \\
0 & -m_1 & * & -m_2 & m_3 & m_4 \\
m_2 & * & m_5 & 0 & * & * \\
* & * & 0 & -m_5 & * & m_6 \\
m_3 & * & * & * & 0 & * \\
m_4 & * & m_6 & * & * & 0
\end{pmatrix}},
\]
where the symbol \(*\) denotes other integral coefficients (not necessarily equal to \(0\) or \(m_1\ldots,m_6\)).
Looking at the possible parities of \(m_1,\ldots,m_6\), we see that \(h = (A \mod 2)\) is conjugate to \(h_2\) in all possible cases.
\endproof

\begin{lemma} \label{lem:order.3.Kloosterman}
An element \(h \in \OO_\hdg^\sharp(T_X)\) of order \(3\) belongs to the conjugacy class of \(h_3\).
\end{lemma}
\proof
Similarly as in \autoref{lem:order.3.Oguiso}, we can assume that \(h = (A\mod 2)\) for some matrix \(A \in \GL_6(\IZ)\) such that \(ATA^\intercal = T\) and \(A^3 = I\).
Again, we let \(U\) and \(W\) be the eigenspace associated with \(1\) and~\(\omega\), respectively. It holds \(\dim(U) + 2\dim(W) = 6\). Moreover, \(U\) is orthogonal to \(W\). 

Investigating all cases with similar computations, we see that it must be \(\dim(U) = \dim(W) = 2\), otherwise \(A\) is forced to have non-integral coefficients.

Looking at the possible parities of the entries of \(A\), we see that \(h = (A \mod 2)\) is always contained in the conjugacy class of~\(h_3\).
\endproof 

\begin{proposition} \label{prop:Hodge.isometries.Kloosterman}
The subgroup \(\OO_\hdg^\sharp(T_X)\) is is a cyclic group of order \(1\), \(2\) or \(3\), 
in which case it is generated by a conjugate of \(h_1\), \(h_2\) or \(h_3\), respectively.
\end{proposition}
\proof
Let \(h\) be a generator of \({\OO_\hdg^\sharp(T_X)}\) and let \(m\) be its order.

By \autoref{prop:HLOY.Appendix.B} it holds \(\varphi(|{\OO_\hdg(T_X)}|) \mid \rank(T_X) = 6\), therefore
\(
    |{\OO_\hdg(T_X)}| \in \set{2,4,6,14,18}.
\)
It follows from \autoref{lem:ker=pm.id.Kloosterman} that \(m = |{\OO_\hdg(T_X)}|/2 \in \set{1,2,3,7,9}\).

We can exclude \(m = 7\) because \(7 \nmid |{\OO(T_X^\sharp)}| = 1440\) and we can exclude \(m = 9\) because we can verify by inspection that \(\OO(T_X^\sharp)\) contains no elements of order~\(9\). 
All in all, \(m \in \set{1,2,3}\) as claimed.

If \(m = 1\), then obviously \(h = h_1\).
If \(m = 2\), then \(h\) is conjugate to \(h_2\) by \autoref{lem:order.2.Kloosterman}. 
If \(m = 3\), then \(h\) is conjugate to \(h_3\) by \autoref{lem:order.3.Kloosterman}. 
\endproof 

\begin{proposition} \label{prop:genus.Kloosterman}
The frame genus \(\cW_X\) contains exactly \(18\) isomorphism classes, listed in \autoref{tab:genus_Kloosterman}, whose Gram matrices are contained in the arXiv ancillary file \verb+genus_Kloosterman.sage+. 
\end{proposition}
\proof 
We compute \(\cW_X\) using Kneser's neighbor method (see \autoref{rem:Kneser.neighbor}). 
One could also have applied the Kneser--Nishiyama method with either \(T_0 = \bA_1^2\bD_4^2\) or \(T_0 = \bA_1^4\bD_6\). (There is a third lattice in the same genus which is not generated by its roots, but this choice makes computations much more difficult, as one cannot use Nishiyama's results \cite{Nishiyama:Jacobian.fibrations.K3.MW}).
The list of lattices found is complete because the mass formula holds (see \autoref{rem:mass}):
\[
    \sum_{i =1}^{18} \frac{1}{|{\OO(W_i)}|} = \frac{1306681}{64210599936000} = \mass(\cW_X). \qedhere
\]
\endproof 

\begin{remark}
As a corollary of \autoref{prop:genus.Kloosterman}, one can easily derive Kloosterman's classification \cite[Thm.~1.1 and Thm.~1.3]{Kloosterman:classification.jac.ell.fibr} of the pairs \((W_\rootlattice,W/W_\rootlattice)\).
Note that \((W_8,W_9)\) and \((W_{10},W_{11})\) correspond to the same pair \((W_\rootlattice,W/W_\rootlattice)\). In Kloosterman's paper these cases are not distinguished. 
In the language of Braun, Kimura and Watari~\cite{Braun.Kimura.Watari:classif.ell.fibr}, in this situation the ``\(\cJ_2(X)\) classification'' is strictly finer than the ``\(\cJ^{\mathrm{(type)}}(X)\) classification''.
For this reason we decided not to follow Kloosterman's numeration.
\end{remark} 

The arXiv ancillary file \verb+groups_K_Kloosterman.gap+ contains the definition of the subgroups 
\[
K_{12}, K_{24}, K_{36}, K_{48}, K_{96}',K_{96}'', K_{120}, K_{144}, K_{240}, K_{720} \subset G.
\]

With the same proof as for \autoref{prop:O^sharp(W).Oguiso} we obtain the following proposition.
\begin{proposition} \label{prop:O^sharp(W).Kloosterman}
The subgroup \(\OO^\sharp(W) \subset G\) is conjugate 
to \(K_{12}\) if \(W = W_9\),  
to \(K_{16}\) if \(W = W_2\),
to \(K_{24}\) if \(W = W_3\),
to \(K_{36}\) if \(W = W_4\),
to \(K_{48}\) if \(W \in \set{ W_5, W_{11} }\),
to \(K_{96}'\) if \(W \in \set{ W_7, W_{13} }\),
to \(K_{96}''\) if \(W \in \set{ W_8, W_{16} }\),
to \(K_{120}\) if \(W = W_{17}\),
to \(K_{144}\) if \(W \in \set{W_1,W_{14},W_{15}}\),
to \(K_{240}\) if \(W = W_{12}\),
to \(K_{720}\) if \(W = W_6\),
and to \(G\) if \(W \in \set{W_{10},W_{18}}\). \qed
\end{proposition}

We now have all ingredients to compute the multiplicities of the frames \(W \in \cW_X\).

\begin{theorem} \label{thm:Kloosterman}
If \(X\) is a K3 surface with transcendental lattice \(T_X \cong \bU(2)^2 \oplus [-2]^2\), then one of the following cases holds:
\[
    |\cJ_X/{\Aut(X)}| = \begin{cases} 482 & \text{if \(|{\OO^\sharp_\hdg(T_X)}| = 1\)}, \\ 274 & \text{if \(|{\OO^\sharp_\hdg(T_X)}| = 2\)}, \\ 172 & \text{if \(|{\OO^\sharp_\hdg(T_X)}| = 3\)}. \end{cases}
\]
\end{theorem}
\proof
We apply \autoref{thm:main_multiplicity} to compute the multiplicities of the frames \(W \in \cW_X\) using either formula~\eqref{eq:Cauchy-Frobenius} or the GAP function \verb+DoubleCosets+.
Thanks to \autoref{prop:Hodge.isometries.Kloosterman} and \autoref{prop:O^sharp(W).Kloosterman}, the multiplicities are given by \(|H\backslash G /K|\), where \(G\) is the group defined in \eqref{eq:group.G.Kloosterman}, \(H\) is the cyclic subgroup generated by \(h_1,h_2\) or \(h_3\) and \(K \in \set{K_{12}, K_{24}, K_{36}, K_{48}, K_{96}',K_{96}'', K_{120}, K_{144}, K_{240}, K_{720}}\).
Our results are listed in \autoref{tab:multiplicities.Kloosterman}.
\endproof

{\small
\begin{table}
    \centering
    \caption{Lattices in the frame genus \(\cW_X\) of a K3 surface \(X\) with transcendental lattice \(T_X \cong \bU(2)^2 \oplus [-2]^2\).}
    \label{tab:genus_Kloosterman}
    \begin{tabular}{llllll}
    \toprule 
    \(W\) & \(W_{\rootlattice}\) & \(W/W_\rootlattice\) & \(|\Delta(W)|\) & \(|{\OO(W)}|\) & \(|{\OO^\sharp(W)}|\) \\
    \midrule 
    \(W_1\) & \(\bD_4^2\bD_6\) & \(0\) & \(108\) & \(122305904640\) & \(144\) \\
    \(W_2\) & \(\bA_1^2\bD_6^2\) & \(0\) & \(124\) & \(33973862400\) & \(16\) \\
    \(W_3\) & \(\bA_1^2\bD_4\bD_8\) & \(0\) & \(140\) &  \(95126814720\) & \(24\) \\
    \(W_4\) & \(\bA_1^3\bD_4\bE_7\) & \(0\) & \(156\) &  \(160526499840\) & \(36\) \\
    \(W_5\) & \(\bA_1^4\bD_{10}\) & \(0\) & \(188\) & \(1426902220800\) & \(48\) \\
    \(W_6\) & \(\bA_1^6\bE_8\) & \(0\) & \(252\) & \(32105299968000\) & \(720\) \\
     \midrule 
    \(W_7\) & \(\bA_1^2\bD_4^3\) & \(\IZ/2\IZ\) & \(76\) & \(2717908992\) & \(96\) \\
    \(W_8\) & \(\bA_1^4\bD_4\bD_6\) & \(\IZ/2\IZ\) & \(92\) & \(6794772480\) & \(96\) \\
    \(W_9\) & \(\bA_1^4\bD_4\bD_6\) & \(\IZ/2\IZ\) & \(92\) & \(849346560\) & \(12\) \\
    \(W_{10}\) & \(\bA_1^6\bD_8\) & \(\IZ/2\IZ\) & \(124\) & \(475634073600\) & \(1440\) \\
    \(W_{11}\) & \(\bA_1^6\bD_8\) & \(\IZ/2\IZ\) & \(124\) & \(15854469120\) & \(48\) \\
    \(W_{12}\) & \(\bA_1^7\bE_7\) & \(\IZ/2\IZ\) & \(140\) & \(89181388800\) & \(240\) \\
    \midrule 
    \(W_{13}\) & \(\bA_1^6\bD_4^2\) & \((\IZ/2\IZ)^2\) & \(60\) & \(226492416\) & \(96\) \\
    \(W_{14}\) & \(\bA_1^8\bD_6\) & \((\IZ/2\IZ)^2\) & \(76\) & \(849346560\) & \(144\) \\
    \midrule 
    \(W_{15}\) & \(\bA_5^2\) & \(\IZ^4\) & \(60\) & \(149299200\) & \(144\) \\
    \(W_{16}\) & \(\bA_3\bA_7\) & \(\IZ^4\) & \(68\) & \(185794560\) & \(96\) \\
    \(W_{17}\) & \(\bA_1\bA_9\) & \(\IZ^4\) & \(92\) & \(1741824000\) & \(120\) \\
    \midrule 
    \(W_{18}\) & \(\bA_3\bE_6\) & \(\IZ^5\) & \(84\) & \(3583180800\) & \(1440\) \\
    \bottomrule 
    \end{tabular}
\end{table}

\begin{table}
    \centering
    \caption{Multiplicities of the frames \(W \in \cW_X\) listed in \autoref{tab:genus_Kloosterman}.}
    \label{tab:multiplicities.Kloosterman}
    \begin{tabular}{c|ccccccccc|c}
    \toprule
        \(|{\OO_\hdg^\sharp(T_X)}|\) & \(W_1\) & \(W_2\) & \(W_3\) & \(W_4\) & \(W_5\) & \(W_6\) & \(W_7\) & \(W_8\) & \(W_9\) & \phantom{\(|{\OO_\hdg^\sharp(T_X)}|\)} \\
    \midrule
    \(1\) & \(10\) & \(90\) & \(60\) & \(40\) & \(30\) & \(2\) & \(15\) & \(15\) & \(120\) & \ldots \\ 
    \(2\) & \(7\) & \(54\) & \(30\) & \(20\) & \(18\) & \(2\) & \(11\) & \(9\) & \(60\) & \ldots \\ 
    \(3\) & \(4\) & \(30\) & \(20\) & \(16\) & \(10\) & \(2\) & \(7\) & \(5\) & \(40\) & \ldots \\
    \bottomrule
    \toprule 
        & \(W_{10}\) & \(W_{11}\) & \(W_{12}\) & \(W_{13}\) & \(W_{14}\) & \(W_{15}\) & \(W_{16}\) & \(W_{17}\) & \(W_{18}\) & sum \\
    \midrule
    \ldots & \(1\) & \(30\) & \(6\) & \(15\) & \(10\) & \(10\) & \(15\) & \(12\) & \(1\) & \(482\) \\
    \ldots & \(1\) & \(18\) & \(3\) & \(11\) & \(7\) & \(7\) & \(9\) & \(6\) & \(1\) & \(274\) \\
    \ldots & \(1\) & \(10\) & \(2\) & \(7\) & \(4\) & \(4\) & \(5\) & \(4\) & \(1\) & \(172\) \\
    \bottomrule 
    \end{tabular}
\end{table}
}

\subsection{Apéry--Fermi pencil} \label{sec:Apery-Fermi}

Let \(X\) be a K3 surface with transcendental lattice
\[
    T_X \cong \bU \oplus [-12].
\]

\begin{remark}
The surface \(X\) is the generic element of a pencil of K3 surfaces studied by Peters and Stienstra~\cite{Peters.Stienstra:pencil.Apery}. 
Elliptic fibrations on this pencil were already classified by Bertin and Lecacheux~\cite{Bertin.Lecacheux:apery-fermi.2-isogenies}, who in particular determined all pairs \((W_\rootlattice, W/W_\rootlattice)\), \(W \in \cW_X\).
\end{remark}

\begin{proposition}
The frame genus \(\cW_X\) contains exactly \(27\) isomorphism classes, listed in \autoref{tab:genus_Apery_Fermi}, whose Gram matrices are contained in the arXiv ancillary file \verb+genus_Apery_Fermi.sage+. 
\end{proposition}
\proof 
The list found by Bertin and Lecacheux~\cite{Bertin.Lecacheux:apery-fermi.2-isogenies} applying the Kneser--Nishiyama method with \(T_0 = \bA_2\bD_5\) is complete because the mass formula holds:
\[
    \sum_{i = 1}^{27} \frac{1}{|{\OO(W_i)}|} = \frac{123970110547}{1110668726060974080000} = \mass(\cW_X). \qedhere 
\]
\endproof

\begin{theorem} \label{thm:Apery-Fermi}
If \(X\) is a K3 surface with transcendental lattice \(T_X \cong \bU \oplus [-12]\), then
\[
    |\cJ_X/{\Aut(X)}| = 32.
\]
\end{theorem}
\proof 
It holds \(|{\OO(T_X^\sharp)}| = 4\). Since \(\rank(T_X)\) is odd, \(|{\fr_X^{-1}(W)}|\) is equal to the index of \(\OO^\sharp(W)\) in \(\OO(T_X^\sharp)\) (\autoref{rmk:rk.TX.odd}).
A finite set of generators of \(\OO(W)\) can be computed using the command \verb+orthogonal_group+ of the Sage class \verb+QuadraticForm+. Their images generate \(\OO^\sharp(W)\).
\endproof 

{\small
\begin{table}[t]
    \centering
    \caption{Lattices in the frame genus \(\cW_X\) of a K3 surface \(X\) with transcendental lattice \(T_X \cong \bU \oplus [-12]\), numbered according to Bertin and Lecacheux (cf. \cite[Table~2 and Table~3]{Bertin.Lecacheux:apery-fermi.2-isogenies}).}
    \label{tab:genus_Apery_Fermi}
    \begin{tabular}{lllllll}
    \toprule 
    \(W\) & \(N_\rootlattice\) & \(W_\rootlattice\) & \(W/W_\rootlattice\) & \(|\Delta(W)|\) & \(|{\OO(W)}|\) & \(|{\fr_X^{-1}(W)}|\) \\
    \midrule
    \(W_3\) & \(\bD_{16}\bE_8\) & \(\bD_{11}\bE_6\)     & \(0\)     & \(292\) & \(8475799191552000\) & \(1\) \\
    \(W_1\) & \(\bE_8^3\) & \(\bA_3\bE_6\bE_8\)         & \(0\)     & \(324\) & \(3467372396544000\) & \(1\) \\
    \midrule
    \(W_7\) & \(\bD_{10}\bE_7^2\) & \(\bA_5\bD_5\bE_7\) & \(\IZ/2\IZ\) & \(196\) & \(16052649984000\) & \(1\) \\
    \midrule
    \(W_{27}\) & \(\bA_7^2\bD_5^2\) & \(\bA_4\bA_7\bD_5\) & \(\IZ\) & \(116\) & \(18579456000\) & \(2\) \\
    \(W_{20}\) & \(\bA_{11}\bD_7\bE_6\) & \(\bA_1^2\bA_2^2\bA_{11}\) & \(\IZ/6\IZ\) & \(148\)
& \(551809843200\) & \(1\) \\
    \midrule
    \(W_{21}\) & \(\bA_{11}\bD_7\bE_6\) & \(\bA_1^2\bA_8\bE_6\) & \(\IZ\) & \(148\) & \(300987187200\) & \(1\) \\
    \(W_{18}\) & \(\bA_{15}\bD_9\) & \(\bA_{12}\bD_4\)  & \(\IZ\) & \(180\) & \(4782351974400\) & \(1\) \\
    \(W_{13}\) & \(\bD_{12}^2\) & \(\bD_9\bD_7\)        & \(\IZ\)   & \(228\) & \(119859786547200\) & \(1\) \\
    \(W_5\) & \(\bD_{16}\bE_8\) & \(\bA_3 \bD_{13}\)    & \(\IZ\)   & \(324\) & \(2448564210892800\) & \(1\) \\
    \(W_6\) & \(\bD_{16}\bE_8\) & \(\bD_8 \bE_8\)       & \(\IZ\)   & \(352\) & \(14383174385664000\) & \(1\) \\
    \(W_2\) & \(\bE_8^3\) & \(\bE_8^2\)                 & \(\IZ\)   & \(480\) & \(1941728542064640000\) & \(2\) \\  
    \(W_{12}\) & \(\bD_{24}\) & \(\bD_{16}\)            & \(\IZ\)   & \(480\) & \(2742391916199936000\) & \(1\) \\
    \midrule 
    \(W_{15}\) & \(\bD_8^3\) & \(\bA_3 \bD_5 \bD_8\)    & \(\IZ \oplus (\IZ/2\IZ)\) & \(164\) & \(951268147200\) & \(1\) \\
    \(W_8\) & \(\bD_{10}\bE_7^2\) & \(\bA_1\bA_5\bD_{10}\) & \(\IZ \oplus (\IZ/2\IZ)\) & \(212\) & \(10701766656000\) & \(1\) \\
    \(W_{16}\) & \(\bD_8^3\) & \(\bD_8^2\)              & \(\IZ \oplus (\IZ/2\IZ)\) & \(224\) & \(106542032486400\) & \(2\) \\
    \(W_9\) & \(\bD_{10}\bE_7^2\) & \(\bA_1^2\bE_7^2\) & \(\IZ \oplus (\IZ/2\IZ)\) & \(256\) & \(269684519731200\) & \(1\) \\
    \(W_{14}\) & \(\bD_{12}^2\) & \(\bD_4 \bD_{12}\)    & \(\IZ \oplus (\IZ/2\IZ)\) & \(288\) & \(753404372582400\) & \(1\) \\
    \(W_4\) & \(\bD_{16}\bE_8\) & \(\bD_{16}\)          & \(\IZ \oplus (\IZ/2\IZ)\) & \(480\) & \(1371195958099968000\) & \(2\) \\
    \midrule
    \(W_{19}\) & \(\bE_6^4\) & \(\bA_2^2\bE_6^2\)       & \(\IZ \oplus (\IZ/3\IZ)\) & \(156\) & \(773967052800\) & \(1\) \\
    \midrule
    \(W_{26}\) & \(\bA_7^2\bD_5^2\) & \(\bA_1^2\bA_7^2\) & \(\IZ \oplus (\IZ/4\IZ)\) & \(116\) & \(52022476800\) & \(1\) \\
    \midrule
    \(W_{25}\) & \(\bA_9^2\bD_6\) & \(\bA_6\bA_9\)      & \(\IZ^2\) & \(132\) & \(73156608000\) & \(1\) \\
    \(W_{22}\) & \(\bA_{11}\bD_7\bE_6\) & \(\bA_8\bD_7\) & \(\IZ^2\) & \(156\) & \(234101145600\) & \(2\) \\
    \(W_{10}\) & \(\bD_{10}\bE_7^2\) & \(\bA_1\bD_7\bE_7\) & \(\IZ^2\) & \(212\) & \(7491236659200\) & \(1\) \\
    \(W_{11}\) & \(\bA_{17}\bE_7\) & \(\bA_1\bA_{14}\)  & \(\IZ^2\) & \(212\) & \(10461394944000\) & \(1\) \\
    \midrule 
    \(W_{24}\) & \(\bD_6^4\) & \(\bA_3\bD_6^2\)         & \(\IZ^2 \oplus (\IZ/2\IZ)\) & \(132\) & \(101921587200\) & \(1\) \\
    \(W_{23}\) & \(\bA_{11}\bD_7\bE_6\) & \(\bA_{11}\bD_4\) & \(\IZ^2 \oplus (\IZ/2\IZ)\) & \(156\) &
\(367873228800\) & \(1\) \\
    \(W_{17}\) & \(\bA_{15}\bD_9\) & \(\bA_{15}\)       & \(\IZ^2 \oplus (\IZ/2\IZ)\) & \(240\) & \(167382319104000\) & \(1\) \\
    \bottomrule
    \end{tabular}
\end{table}
}

\bibliographystyle{amsplain}
\bibliography{references}

\end{document}